\documentclass[11pt,a4paper]{article}
\usepackage[a4paper,left=1.25in,right=1.25in,top=1.1in,bottom=1.1in]{geometry}
\usepackage{amsmath,amssymb,amsthm}
\usepackage{booktabs,graphicx,array}
\usepackage{algorithm}
\usepackage{algpseudocode}
\usepackage[dvipsnames]{xcolor}
\usepackage[colorlinks=true,linkcolor=blue!50!black,citecolor=blue!50!black,urlcolor=blue!50!black]{hyperref}
\usepackage{enumitem}

\newtheorem{theorem}{Theorem}[section]
\newtheorem{proposition}[theorem]{Proposition}
\newtheorem{corollary}[theorem]{Corollary}
\newtheorem{lemma}[theorem]{Lemma}
\theoremstyle{remark}
\newtheorem{remark}[theorem]{Remark}
\DeclareMathOperator*{\argmin}{arg\,min}
\newcommand{\R}{\mathbb{R}}

\newcommand{\norm}[1]{\lVert #1\rVert}
\newcommand{\ip}[2]{\langle #1,#2\rangle}
\title{\textbf{SPIKE-Accelerated Split-Bregman/ADMM Methods for Obstacle Problems and American Option Pricing under Jump-Diffusion Models}}
\author{}
\date{}

\begin{document}
	\maketitle
	\vspace{-2.5em}
	\begin{center}
		\large
		Amboru Yalamanda$^1$ and Peeyush Singh$^2$ \\[6pt]
		\small
		$^1$ Department of Mathematics and Statistics, VIT-AP UNIVERSITY, Amaravathi, Andhra Pradesh, India\\
		amboru.24phd7073@vitap.ac.in\\[4pt]
		$^2$ Department of Mathematics and Statistics, VIT-AP UNIVERSITY, Amaravathi, Andhra Pradesh, India\\
		peeyush.singh@vitap.ac.in
	\end{center}
	\vspace{1em}
	
	\begin{abstract}
		\noindent This paper presents a parallel split Bregman/alternating direction method of multipliers (ADMM) \cite{goldstein,tran} for obstacle problems and for the linear complementarity problems (LCPs) that arise in the pricing of American options under the Black--Scholes and Merton jump-diffusion \cite{merton} models. The auxiliary variable of ADMM decouples the inequality constraint from the differential operator, so that every iteration consists of a linear solve with a fixed banded matrix and a pointwise projection. The matrix is factorized once by the SPIKE algorithm \cite{spike}, as in the parallel direct approach of \cite{singh}, and reused in all iterations and time steps in a Fortran~90 MPI implementation. Convergence is proved by a Lyapunov functional for every positive augmentation parameter, for the symmetric obstacle problem and the non-symmetric option LCP alike, and the iterates do not depend on the number of processes. The method therefore needs a single factorization, no tuned relaxation parameter and only pointwise, perfectly parallel updates. One- and two-dimensional obstacle problems, the Black--Scholes call and American put, and European and American options under Merton's model are computed and compared with exact, closed-form, binomial and published values. In two dimension a speed-up of $3.92$ on four processes is obtained. It is also shown that the augmentation parameter has to be scaled to the diagonal of the time-stepping matrix when only a few inner iterations are used.
	\end{abstract}
	
	
	%======================================================================
	\section{Introduction}
	Obstacle problems describe a membrane that is held above a rigid obstacle and occur in contact mechanics, flow through porous media and free-boundary problems \cite{tran,kinderlehrer}. A problem of the same structure arises in finance. Black and Scholes \cite{bs} formulated a diffusion model for the price of an option. Because an American option can be exercised at any time before maturity \cite{merton73}, its price is bounded below by the payoff. The price is then the solution of a variational inequality, in which the payoff acts as an obstacle \cite{ikonen,maurya}. In the Black--Scholes model the volatility is constant, which does not reproduce the volatility smiles and heavy tails that are observed in the markets. The jump-diffusion models of Merton \cite{merton} and Kou \cite{kou}, in which the log-price is driven by a Brownian motion and a Poisson process with normally or double-exponentially distributed jump sizes, address this, but the integral term that describes the jumps makes the discrete operator dense. After finite-difference discretization, obstacle problems and American options lead to a convex problem or a linear complementarity problem (LCP) with a banded matrix, tridiagonal in one dimension and block tridiagonal with half-bandwidth $n_x$ for the five-point stencil in two dimensions.
	
	Many numerical methods have been proposed for the jump-diffusion PIDE. Cont and Voltchkova \cite{cont} treated the integral term explicitly and the differential part implicitly. D'Halluin et al.\ \cite{dhalluin} solved the dense system of the implicit scheme by a fixed-point iteration, and Kwon and Lee \cite{kwon} developed a second-order tridiagonal method for American options. For the LCP, projected successive over-relaxation (PSOR) \cite{cryer} is the classical solver; it is sequential and needs a relaxation parameter that has to be tuned. Ikonen and Toivanen \cite{ikonen} introduced operator splitting, in which the Black--Scholes operator is separated from the constraint so that each time step needs one tridiagonal solve and a projection, and Maurya et al.\ \cite{maurya} combined this idea with implicit--explicit Runge--Kutta methods for the Merton and Kou models. Penalty methods replace the constraint by a penalty term whose solution depends on the penalty parameter, and active-set methods need few iterations but change the matrix whenever the contact set changes, so that the matrix has to be refactorized. Tran, Schaeffer, Feldman and Osher \cite{tran} replaced the constraint by an exact $L^1$ penalty and solved the obstacle problem by the split Bregman method \cite{goldstein}, which is equivalent to the alternating direction method of multipliers (ADMM) \cite{glowinski,gabay,boyd}; it is connected to Bregman iteration \cite{yin} and to Douglas--Rachford splitting \cite{douglas,lions,eckstein}.
	
	The linear solve is the dominant cost of all of these methods. SPIKE \cite{spike,polizzi07} is a divide-and-conquer solver for banded systems that is related to the partitioning ideas of Sameh and Kuck \cite{sameh78}: the matrix is cut into diagonal blocks that are factorized independently, and the coupling between the blocks is condensed into a small reduced system. SPIKE has been used for obstacle problems together with a projection step by Singh and Yalamanda \cite{singh}, who proposed parallel direct and iterative methods for general obstacle and free boundary problems, including a SPIKE-based direct solver in one dimension and a directional splitting approach in two and three dimensions. This paper combines ADMM and SPIKE. In the ADMM formulation an auxiliary variable decouples the constraint from the operator, so that every iteration consists of one linear solve with a \emph{fixed} banded matrix and a closed-form pointwise update. The matrix is constant over the iterations and, for options, over the time steps, so that the SPIKE factorization, the spikes and the reduced system are computed once and reused in every iteration of a distributed-memory (MPI) implementation. To our knowledge, the combination of a fixed-matrix splitting with a parallel banded solver has not been used for these problems.
	
	The approach has the following advantages, which are made precise in Sections~\ref{sec:why}--\ref{sec:results}. A single factorization serves all iterations and time steps, and nothing is refactorized when the contact set changes. The iteration converges for every positive augmentation parameter, so that, in contrast to PSOR, no relaxation parameter has to be tuned. The projection and the multiplier update are pointwise and perfectly parallel. The same algorithm and the same convergence proof cover the symmetric obstacle problem and the non-symmetric LCP of an American option. Finally, since SPIKE is an exact solver, the iterates do not depend on the number of processes. The convergence proof of Section~\ref{sec:conv} is self-contained and is complemented by the exact-penalty result of \cite{tran}. The method is tested on one- and two-dimensional obstacle problems, on the Black--Scholes call and American put, and on European and American options under Merton's model, including the analytic tails of the truncated jump integral. A practical finding is that, if only a fixed number of inner iterations is performed in every time step, the augmentation parameter has to be scaled to the diagonal of the time-stepping matrix; otherwise the price is systematically too low.
	
	The rest of the paper is organized in the following manner. Section~\ref{sec:sb} introduces the split Bregman and ADMM methods and Section~\ref{sec:model} the mathematical models. The method for obstacle problems is formulated in Section~\ref{sec:method} and its use for option pricing in Section~\ref{sec:option}. Section~\ref{sec:why} discusses the advantages of the approach, and Sections~\ref{sec:theory} and~\ref{sec:conv} contain the mathematical theory and the convergence proof. The SPIKE algorithm and its parallel implementation are described in Section~\ref{sec:spike}, and the parallel performance is reported in Section~\ref{sec:perf}. The test problems and the numerical results are given in Sections~\ref{sec:problems} and~\ref{sec:results}, and a summary and the conclusions are given in Section~\ref{sec:concl}.
	
	%======================================================================
	\section{Split Bregman and ADMM methods}\label{sec:sb}
	The linear complementarity problems of the Introduction are convex problems with a simple constraint set. This section recalls the splitting that separates the constraint from the differential operator; it is the basis of the method of Section~\ref{sec:method} and of the convergence proof of Section~\ref{sec:conv}.
	Consider the problem
	\begin{equation}\label{eq:gen}
		\min_{u,v}\ f(u)+g(v)\quad\text{subject to}\quad u=v,
	\end{equation}
	with $f$ smooth and strongly convex and $g$ closed, convex and possibly extended-valued. ADMM \cite{glowinski,gabay,boyd} minimizes the augmented Lagrangian alternately in $u$ and $v$ and then updates the multiplier. With the scaled multiplier $d=y/\lambda$ and penalty parameter $\lambda>0$, one iteration is
	\begin{align}
		u^{k+1}&=\argmin_u\ f(u)+\tfrac\lambda2\norm{u-v^k+d^k}_2^2, \label{eq:admm1}\\
		v^{k+1}&=\argmin_v\ g(v)+\tfrac\lambda2\norm{u^{k+1}+d^k-v}_2^2, \label{eq:admm2}\\
		d^{k+1}&=d^k+u^{k+1}-v^{k+1}. \label{eq:admm3}
	\end{align}
	The split Bregman method of Goldstein and Osher \cite{goldstein} is the same iteration written with Bregman variables: the quantity $d^k$ is the accumulated constraint residual that is added back to the data of the quadratic term, as in Bregman iteration \cite{yin}. For the obstacle problem $f(u)=\tfrac12 u^{\top}Au-F^{\top}u$ is quadratic, so \eqref{eq:admm1} is a linear system with the constant matrix $A+\lambda I$, and if $g$ is the indicator function of $C=\{v\ge\varphi\}$, then \eqref{eq:admm2} is the projection
	\begin{equation}\label{eq:proj}
		v^{k+1}=\max\bigl(\varphi,\ u^{k+1}+d^k\bigr).
	\end{equation}
	The $L^1$ penalty $\mu(\varphi-v)_+$ of \cite{tran} gives instead a shrinkage; their Theorem~3.1 states that the two formulations have the same minimizer when $\mu$ is large enough (we use only the hard-constraint form, in which $\mu$ does not enter the iteration). The projection \eqref{eq:proj} is the one-sided shrinkage of this formulation with an infinite penalty. In this sense the exact-penalty bound of \cite{tran} and the convergence theorem of Section~\ref{sec:conv} complement each other: the iteration converges, for every $\lambda>0$, to the solution of the constrained problem, which by \cite{tran} is also the solution of the penalized problem when $\mu$ is large enough.
	
	%======================================================================
	\section{Mathematical models}\label{sec:model}
	Having fixed the splitting, we now state the two problem classes to which it is applied: the continuous obstacle problem and the LCP of an American option under the Black--Scholes and Merton models. Their discretizations give the matrices that the splitting and SPIKE act on in Sections~\ref{sec:method} and~\ref{sec:option}.
	\subsection{The obstacle problem}
	Let $\Omega\subset\R^d$ be bounded and $\varphi$ a given obstacle. We seek
	\begin{equation}\label{eq:obst}
		u=\argmin_{v\ge\varphi}\ \tfrac12\int_\Omega|\nabla v|^2\,dx,\qquad v=0\ \text{on }\partial\Omega,
	\end{equation}
	equivalently $\int_\Omega\nabla u\cdot\nabla(v-u)\,dx\ge0$ for all $v\ge\varphi$. The obstacles of \cite{tran}, Section~6.1, are, on $[0,1]$ and extended symmetrically from $[0,0.5]$,
	\begin{align}
		\varphi_1(x)&=\begin{cases}100x^2,&0\le x\le0.25,\\ 100x(1-x)-12.5,&0.25\le x\le0.5,\end{cases}\label{eq:phi1}\\
		\varphi_2(x)&=\begin{cases}10\sin(2\pi x),&0\le x\le0.25,\\ 5\cos(\pi(4x-1))+5,&0.25\le x\le0.5,\end{cases}\label{eq:phi2}
	\end{align}
	with exact solutions
	\begin{equation}\label{eq:exact}
		u_1(x)=\begin{cases}(100-50\sqrt2)\,x,&0\le x\le\frac1{2\sqrt2},\\ 100x(1-x)-12.5,&\frac1{2\sqrt2}\le x\le0.5,\end{cases}\qquad
		u_2(x)=\begin{cases}10\sin(2\pi x),&0\le x\le0.25,\\ 10,&0.25\le x\le0.5.\end{cases}
	\end{equation}
	In two dimensions we use the non-smooth obstacle $\varphi_3$ on $[0,1]^2$ of \cite{tran}, Eq.~(6.4): $\varphi_3=5.0$ on the diamond $|x-0.6|+|y-0.6|<0.04$, $4.5$ on the disc $(x-0.6)^2+(y-0.25)^2<0.001$, $4.5$ on the segment $y=0.57$, $0.075<x<0.13$, and $0$ elsewhere.
	
	\subsection{American options under Black--Scholes and jump-diffusion models}
	Let $S$ be the price of the underlying asset, $E$ the strike price, $T$ the maturity, $r$ the risk-free interest rate and $\sigma$ the volatility. With $\tau=T-t$, the value $u(\tau,\cdot)$ of an American option with payoff $g$ satisfies the complementarity system \cite{ikonen,maurya}
	\begin{equation}\label{eq:lcp}
		u_\tau-\mathcal L u\ge0,\qquad u\ge g,\qquad (u_\tau-\mathcal L u)(u-g)=0 .
	\end{equation}
	For the Black--Scholes model \cite{bs}, $\mathcal Lu=\tfrac12\sigma^2S^2u_{SS}+rSu_S-ru$. For jump-diffusion models, in the log-price $x=\ln(S/E)$, $\mathcal L=\mathcal D+\lambda_J\mathcal I$ with
	\begin{align}
		\mathcal Du&=\tfrac12\sigma^2u_{xx}+\bigl(r-\tfrac{\sigma^2}{2}-\lambda_J\kappa\bigr)u_x-(r+\lambda_J)u,\\
		\mathcal Iu(x)&=\int_{-\infty}^{\infty}u(y)\,\gamma(y-x)\,dy,
	\end{align}
	where $\lambda_J$ is the jump intensity (the mean number of jumps per unit time), $\gamma$ is the density of the jump size and $\kappa=\int(e^{y}-1)\gamma(y)\,dy$ \cite{maurya}. In Merton's model \cite{merton} $\gamma$ is the $N(\mu_J,\sigma_J^2)$ density, so $\kappa=e^{\mu_J+\sigma_J^2/2}-1$. The drift contains $\lambda_J\kappa$; omitting it changes the model.
	
	On a truncated domain $(x_{\min},x_{\max})$ the integral is split as $\mathcal Iu=\int_{x_{\min}}^{x_{\max}}u\,\gamma(y-x)\,dy+\mathcal R^*(x,\tau)$, where $\mathcal R^*$ integrates the asymptotic option value outside the domain \cite{maurya,cont,dhalluin}. With $\Phi$ the standard normal distribution function, a put has the asymptotic value $Ee^{-r\tau}-S$ (European) or $E-S$ (American) at the left end, and \cite{maurya}
	\begin{align}
		\mathcal R^*_{\rm Eur}(x,\tau)&=Ee^{-r\tau}\,\Phi\!\Bigl(\tfrac{x_{\min}-x-\mu_J}{\sigma_J}\Bigr)-Ee^{x+\mu_J+\sigma_J^2/2}\,\Phi\!\Bigl(\tfrac{x_{\min}-x-\mu_J-\sigma_J^2}{\sigma_J}\Bigr),\label{eq:tailE}\\
		\mathcal R^*_{\rm Am}(x)&=E\,\Phi\!\Bigl(\tfrac{x_{\min}-x-\mu_J}{\sigma_J}\Bigr)-Ee^{x+\mu_J+\sigma_J^2/2}\,\Phi\!\Bigl(\tfrac{x_{\min}-x-\mu_J-\sigma_J^2}{\sigma_J}\Bigr).\label{eq:tailA}
	\end{align}
	For a call the same derivation on the right tail, with the asymptotic value $Ee^{x}-Ee^{-r\tau}$ at $x_{\max}$, gives
	\begin{equation}\label{eq:tailC}
		\mathcal R^*_{\rm call}(x,\tau)=Ee^{x+\mu_J+\sigma_J^2/2}\Bigl[1-\Phi\!\Bigl(\tfrac{x_{\max}-x-\mu_J-\sigma_J^2}{\sigma_J}\Bigr)\Bigr]-Ee^{-r\tau}\Bigl[1-\Phi\!\Bigl(\tfrac{x_{\max}-x-\mu_J}{\sigma_J}\Bigr)\Bigr].
	\end{equation}
	Equation~\eqref{eq:tailC} is our own analogue of \eqref{eq:tailE}; its use is checked against published call prices in Section~\ref{sec:results}. For Kou's double-exponential model the corresponding tail has the same structure with an exponential jump density; this model is not treated numerically here.
	
	\subsection{Notation}\label{sec:notation}
	The symbols used in the paper are collected in Table~\ref{tab:notation}. A few letters have different meanings in different parts of the paper; the meaning is then fixed by the section named in the table. In particular, $\lambda$ is the ADMM penalty parameter in the obstacle problems and $\lambda_J$ the jump intensity in the option problems (the penalty parameter of the option problems is $\rho$).
	\begin{table}[h]
		\centering\small
		\caption{Notation.}\label{tab:notation}
		\renewcommand{\arraystretch}{1.15}
		\begin{tabular}{@{}p{0.2\linewidth}p{0.74\linewidth}@{}}
			\toprule
			Symbol & Meaning\\
			\midrule
			$\Omega,\ \varphi,\ u$ & domain, obstacle and solution of the obstacle problem (Section~\ref{sec:model})\\
			$h,\ n,\ n_x,\ N$ & mesh width; number of interior grid points in 1D; per direction in 2D; number of unknowns, $N=n$ or $n_x^2$\\
			$K_h,\ A$ & discrete Laplacian (symmetric positive definite); a generic symmetric positive definite matrix\\
			$M=I+\Delta t\,A_h$ & non-symmetric time-stepping matrix of the option problem (Section~\ref{sec:option}); $A_h$ discretizes $-\mathcal D$\\
			$f,\ g,\ C$ & in Section~\ref{sec:sb}: smooth and non-smooth part of the splitting and the constraint set $C=\{v\ge\varphi\}$; elsewhere $g$ is the option payoff\\
			$u,\ v,\ d$ & primal variable, auxiliary (split) variable and scaled multiplier of ADMM\\
			$\lambda,\ \rho$ & ADMM penalty (augmentation) parameter: $\lambda$ for obstacle problems, $\rho$ for option problems\\
			$k,\ \mathrm{tol},\ k_{\max}$ & iteration counter, stopping tolerance and iteration limit\\
			$S,\ E,\ T,\ \tau=T-t$ & asset price, strike, maturity, time to maturity\\
			$r,\ \sigma$ & risk-free rate and volatility\\
			$x=\ln(S/E)$ & log-price\\
			$\mathcal L,\ \mathcal D,\ \mathcal I$ & full operator, differential part and jump integral of the option problem\\
			$\lambda_J,\ \gamma,\ \kappa$ & jump intensity, density of the jump size, compensator $\kappa=\int(e^y-1)\gamma(y)\,dy$\\
			$\mu_J,\ \sigma_J$ & mean and standard deviation of the log-jump size in Merton's model\\
			$\Phi,\ \mathcal R^*$ & standard normal distribution function; tail of the jump integral outside the truncated domain\\
			$\Delta t,\ \alpha,\ \beta$ & time step; $\alpha=\sigma^2/2$, $\beta=r-\sigma^2/2-\lambda_J\kappa$\\
			$M,\ N$ in Sections~\ref{sec:problems}, \ref{sec:results} & number of space intervals and of time steps of the option grid (not the matrix $M$)\\
			$p,\ P$ & number of slabs/MPI processes ($P$ in the timing tables)\\
			$B=DS,\ B_j,\ V_j,\ W_j$ & banded matrix $K_h+\lambda I$ or $M+\rho I$, its diagonal blocks, and the left and right spikes (Section~\ref{sec:spike}); $S$ here is the spike matrix, not the asset price\\
			$m,\ b$ & width of the coupling blocks; half-bandwidth of $B$\\
			$\norm{\cdot},\ \ip{\cdot}{\cdot}$ & Euclidean norm and inner product; $\norm{\cdot}_\infty$ maximum norm\\
			$m,\ L,\ E(u),\ N_C$ & in Sections~\ref{sec:theory}, \ref{sec:conv}: $m=\lambda_{\min}(K)$ (modulus of strong monotonicity), $L=\lambda_{\max}(K)$, energy $E(u)=\frac12u^{\top}Ku$ (not the strike), normal cone of $C$\\
			$K,\ k$ & in the option problems $K=E$ is the strike; $K$ is also the matrix of Section~\ref{sec:conv}; $k$ is the iteration index\\
			$\zeta_1,\zeta_2,\ p$ & parameters of Kou's jump density (Problem 7)\\
			\bottomrule
		\end{tabular}
	\end{table}
	
	%======================================================================
	\section{ADMM--SPIKE method for obstacle problems}\label{sec:method}
	This section applies the splitting of Section~\ref{sec:sb} to the discrete obstacle problem of Section~\ref{sec:model} and shows how SPIKE enters the linear solve. The same construction is extended to option pricing in Section~\ref{sec:option}.
	Take the uniform grid $x_i=ih$, $h=1/(n+1)$, in 1D, or the $n_x\times n_x$ interior grid with $h=1/(n_x+1)$ in 2D with lexicographic ordering. Let $K_h=-\Delta_h$ be the standard three-point or five-point matrix scaled by $h^{-2}$ (so that $K_h$ approximates $-\Delta$). It is symmetric positive definite, an $M$-matrix \cite{berman}, tridiagonal in 1D, and block tridiagonal with half-bandwidth $n_x$ in 2D. The discrete problem is $\min_{u\ge\varphi}\frac12u^{\top}K_hu$. The scaled ADMM iteration \eqref{eq:admm1}--\eqref{eq:admm3} is:
	\begin{algorithm}[t]
		\caption{ADMM-SPIKE for the discrete obstacle problem}\label{alg:obst}
		\begin{algorithmic}[1]
			\Require $K_h$, obstacle $\varphi$, penalty $\lambda>0$, number of processes $p$, tolerance $\mathrm{tol}$, iteration limit $k_{\max}$
			\Ensure approximate solution $v$ of $\min_{u\ge\varphi}\tfrac12u^{\top}K_hu$
			\State \textbf{Setup (once, in parallel):} partition the grid rows into $p$ slabs and $B=K_h+\lambda I$ into diagonal blocks $B_1,\dots,B_p$ (rank $j$ owns $B_j$); each rank factorizes $B_j$ and computes its spikes $V_j,W_j$ \eqref{eq:spikes}; rank~0 gathers the interface blocks and factorizes the reduced matrix of \eqref{eq:red1}--\eqref{eq:red2}
			\State $v^0\leftarrow\varphi$, \ $d^0\leftarrow0$, \ $k\leftarrow0$
			\Repeat
			\State $u^{k+1}\leftarrow B^{-1}\,\lambda\,(v^k-d^k)$ \Comment{SPIKE solve, Section~\ref{sec:spike}: steps 1--4 of the MPI structure}
			\State $v^{k+1}\leftarrow\max(\varphi,\ u^{k+1}+d^k)$ \Comment{projection onto $\{v\ge\varphi\}$, componentwise}
			\State $d^{k+1}\leftarrow d^k+u^{k+1}-v^{k+1}$ \Comment{multiplier update}
			\State $k\leftarrow k+1$
			\Until{$\norm{u^{k}-u^{k-1}}_\infty<\mathrm{tol}$ \textbf{or} $k=k_{\max}$}
			\State \Return $v^{k}$ \Comment{feasible, $v^k\ge\varphi$}
		\end{algorithmic}
	\end{algorithm}
	
	\noindent The stopping test is that of \cite{tran}; Section~\ref{sec:results} shows that it does not control the error in 2D. Replacing the $L^1$ penalty of \cite{tran} by the indicator of $\{u\ge\varphi\}$ turns their shrink operator into the projection in step~5. The matrix $K_h+\lambda I$ is symmetric positive definite and strictly diagonally dominant, and it does not change from one iteration to the next. This constancy is what makes SPIKE attractive: its factorization, which is the expensive part, is paid once, and every iteration only needs triangular solves and one small reduced solve.
	
	%======================================================================
	\section{ADMM--SPIKE for American option pricing under jump-diffusion models}\label{sec:option}
	The obstacle-problem algorithm of Section~\ref{sec:method} carries over to the option LCP of Section~\ref{sec:model} once the time stepping is fixed. The new points are that the matrix is non-symmetric, that a time loop surrounds the ADMM iteration, and that the augmentation parameter has to match the scaling of the time-stepping matrix.
	\subsection{Time discretization and the linear complementarity problem}
	We use implicit Euler in $\tau$ with step $\Delta t$ and treat the jump integral explicitly, $\mathcal Iu$ being evaluated at the previous time level; this is the implicit--explicit idea of \cite{maurya} without Runge--Kutta stages. In the log-price, with $\alpha=\sigma^2/2$ and $\beta=r-\sigma^2/2-\lambda_J\kappa$, central differences give the tridiagonal rows
	\[
	\Bigl(-\Delta t\bigl(\tfrac{\alpha}{h^2}-\tfrac{\beta}{2h}\bigr),\ 1+\Delta t\bigl(\tfrac{2\alpha}{h^2}+r+\lambda_J\bigr),\ -\Delta t\bigl(\tfrac{\alpha}{h^2}+\tfrac{\beta}{2h}\bigr)\Bigr).
	\]
	These rows form the matrix $M=I+\Delta t\,A_h$, where $A_h$ is the finite-difference matrix of $-\mathcal D$ with the zeroth-order term $r+\lambda_J$. In the Black--Scholes price variable $S_j=j\,\Delta S$ the rows are $\bigl(-\tfrac12\Delta t(\sigma^2j^2-rj),\ 1+\Delta t(\sigma^2j^2+r),\ -\tfrac12\Delta t(\sigma^2j^2+rj)\bigr)$. Each time step is the LCP
	\begin{equation}\label{eq:dlcp}
		\resizebox{\linewidth}{!}{$\displaystyle u^{n+1}\ge g,\quad Mu^{n+1}-b^n\ge0,\quad (Mu^{n+1}-b^n)^{\top}(u^{n+1}-g)=0,\qquad b^n=u^n+\Delta t\,\lambda_J\,\mathcal I_h(u^n)+(\text{boundary terms}),$}
	\end{equation}
	where $\mathcal I_h$ is the trapezoidal approximation of the integral over the domain plus the tail $\mathcal R^*$. (Without a jump term, $\lambda_J=0$, $b^n=u^n$.) The matrix $M$ is tridiagonal and not symmetric.
	
	\subsection{ADMM iteration for one time step}
	\begin{algorithm}[t]
		\caption{ADMM-SPIKE for the time stepping of the LCP \eqref{eq:dlcp}}\label{alg:opt}
		\begin{algorithmic}[1]
			\Require $M$, payoff $g$, penalty $\rho>0$, time step $\Delta t$, number of steps $N$, number of processes $p$, tolerance $\mathrm{tol}$, sweep limit $k_{\max}$
			\Ensure option values $u^{N}$ at $\tau=N\Delta t$
			\State \textbf{Setup (once for the whole time loop):} SPIKE factorization of $B=M+\rho I$ as in Algorithm~\ref{alg:obst}
			\State $u^0\leftarrow g$, \ $v\leftarrow g$, \ $d\leftarrow0$
			\For{$n=0,1,\dots,N-1$}
			\State form $b^n$ from $u^n$ \Comment{jump integral and tails at $\tau_n$, boundary values at $\tau_{n+1}$}
			\State $k\leftarrow0$, \ $v^0\leftarrow v$, \ $d^0\leftarrow d$ \Comment{warm start from the previous step}
			\Repeat
			\State $u^{k+1}\leftarrow B^{-1}\bigl(b^n+\rho\,(v^k-d^k)\bigr)$ \Comment{SPIKE solve}
			\State $v^{k+1}\leftarrow\max(g,\ u^{k+1}+d^k)$
			\State $d^{k+1}\leftarrow d^k+u^{k+1}-v^{k+1}$
			\State $k\leftarrow k+1$
			\Until{$\bigl(\norm{v^{k}-v^{k-1}}_\infty<\mathrm{tol}$ \textbf{and} $\norm{u^{k}-v^{k}}_\infty<\mathrm{tol}\bigr)$ \textbf{or} $k=k_{\max}$}
			\State $u^{n+1}\leftarrow v^{k}$, \ $v\leftarrow v^{k}$, \ $d\leftarrow d^{k}$
			\EndFor
			\State \Return $u^{N}$
		\end{algorithmic}
	\end{algorithm}
	Theorem~\ref{thm:admm} applies to each time step with the warm-start values as starting data (it holds from $k=1$ on for arbitrary $v^0,d^0$). The inner iteration must be run to convergence in every time step. If it is stopped after a small fixed number of sweeps with a large $\rho$, the contraction factor per sweep is close to $\rho/(1+\rho)$ and only a fraction of each time step is applied, so the option value is systematically too low. This is demonstrated in Section~\ref{sec:results}. Note also that the stopping test must contain the change $\norm{v^{k+1}-v^k}$: the residual $\norm{u^{k+1}-v^{k+1}}$ alone vanishes after one sweep whenever the constraint is inactive and then does not measure convergence.
	
	\subsection{Choice of the augmentation parameter}
	The matrix $M=I+\Delta t\,A_h$ has a diagonal that is not constant: on the Black--Scholes grid of Section~\ref{sec:problems} its entries $1+\Delta t(\sigma^2j^2+r)$ range from $1.0$ to $11.4$ (median $3.6$). In the $u$-update the penalty term acts on node $j$ with the weight $\rho/(m_{jj}+\rho)$, so the contraction per sweep is of order $\rho/(m_{jj}+\rho)$. The augmentation parameter must therefore be of the size of the diagonal of $M$, i.e.\ $\rho=O(1)$ in the scaling used here. A value much larger than the diagonal, such as $\rho=50$, makes every sweep a small correction of order $1/\rho$ and a fixed per-step iteration budget is then far too small. Table~\ref{tab:rho} in Section~\ref{sec:p5} quantifies this.
	
	\subsection{Constancy of the system matrix}
	In Algorithm~\ref{alg:opt} the matrix $M+\rho I$ depends on $\Delta t$, the grid and the model parameters only, so a single SPIKE setup serves all time steps and all inner iterations; the jump integral enters only through the right-hand side. For the jump integral, $\mathcal I_h u$ is a dense matrix--vector product, $O(M^2)$ per step (or $O(M\log M)$ with an FFT since the kernel is a convolution in the log-price \cite{dhalluin}); it is evaluated on rank~0 or distributed by rows.
	
	%======================================================================
	\section{Rationale and properties of the approach}\label{sec:why}
	This section motivates the combination of the ADMM splitting with a SPIKE-based direct solver, compares it with the classical alternatives for the problems of Section~\ref{sec:model}, and summarizes its properties.
	
	\subsection{Split Bregman/ADMM as the constraint-handling scheme}
	The constraint $u\ge\varphi$ makes the discrete problem nonlinear and non-smooth. The ADMM splitting of Section~\ref{sec:sb} removes this difficulty by separating the problem into two simple subproblems: a \emph{linear} system with a constant, symmetric positive definite matrix (the $u$-update \eqref{eq:it1}), and a \emph{pointwise} projection onto $\{v\ge\varphi\}$ (the $v$-update \eqref{eq:it2}). Consequently,
	\begin{itemize}[itemsep=1pt]
		\item every iteration costs one linear solve with a fixed matrix plus $O(N)$ work for the projection and the multiplier update, and the projection is perfectly parallel;
		\item the iteration converges for \emph{every} penalty parameter $\lambda>0$ (Theorem~\ref{thm:admm}), whereas projected successive over-relaxation (PSOR) \cite{cryer} and Gauss--Seidel type sweeps are sequential by nature and require a relaxation parameter $\omega\in(1,2)$ that has to be tuned for each problem;
		\item the same iteration, with the same convergence proof, covers the symmetric obstacle problem and the non-symmetric LCP of an American option (Section~\ref{sec:conv}), and does not require a penalty reformulation of the constraint, in contrast to penalty methods whose solution depends on the penalty parameter \cite{tran}.
	\end{itemize}
	
	\subsection{SPIKE as the parallel banded solver}
	The linear system of the $u$-update is banded, and its matrix does not change from one iteration (or time step) to the next. A serial banded factorization of the 2D problem, with $N=n_x^2$ unknowns and half-bandwidth $n_x$, costs $O(n_x^4)$ operations \cite{golub} and is the dominant cost of the method. SPIKE \cite{spike,polizzi07} distributes this work over $p$ processes: the diagonal blocks are factorized independently, the coupling is condensed into a small reduced system of dimension $2m(p-1)$, and, because the matrix is constant, factorization, spikes and reduced system are computed only once and reused in every iteration (Section~\ref{sec:spike}). The per-iteration cost is then limited to local triangular solves, one reduced solve and the exchange of the interface values.
	
	In one space dimension the matrix is tridiagonal, a serial Thomas solve costs $O(n)$ and cannot be beaten by SPIKE at the sizes considered here. The one-dimensional problems therefore serve to verify the algorithm and its parallel structure; no speed-up over the Thomas algorithm is claimed (Section~\ref{sec:perf}).
	
	\subsection{Exact versus inexact inner solvers}
	Because SPIKE with the full spikes and the full reduced system is an \emph{exact} solver, the iterates do not depend on the number of processes or on the partition (Proposition~\ref{prop:spike}), and the convergence theory of ADMM carries over to the parallel code without change. This is not true for iterative or domain-decomposition solvers that solve the linear system only approximately, for which the inexactness enters the ADMM iteration and has to be controlled separately.
	
	\subsection{Comparison with alternative methods}
	Table~\ref{tab:why} summarizes the properties of the methods mentioned above for the problems of this paper.
	\begin{table}[t]
		\centering\small
		\caption{Qualitative comparison of solution methods for the discrete obstacle problem and the option LCP.}\label{tab:why}
		\begin{tabular}{@{}p{2.9cm}p{3.1cm}p{3.5cm}p{3.5cm}@{}}
			\toprule
			Method & Tuning & Parallel structure & Proven convergence\\
			\midrule
			PSOR \cite{cryer} & $\omega\in(1,2)$ & sequential sweeps & yes (SPD matrices)\\
			Active-set / policy iteration & none & one banded solve per outer iteration, varying matrix & finite termination for $M$-matrices\\
			$L^1$-penalty split Bregman \cite{tran} & penalty parameter, shrinkage & not addressed there & yes for the penalized problem\\
			ADMM + SPIKE (this paper) & $\lambda>0$ (any value) & constant matrix, setup once, exact parallel solve & yes, for every $\lambda>0$ (Theorem~\ref{thm:admm})\\
			\bottomrule
		\end{tabular}
	\end{table}
	
	\subsection{Principal properties}
	Collecting the points above, the approach has the following properties for obstacle problems and for the LCPs of option pricing.
	\begin{enumerate}[itemsep=1pt]
		\item \emph{A fixed matrix.} The matrix $K_h+\lambda I$ (or $M+\rho I$) is the same in all iterations and, for options, in all time steps; the SPIKE factorization, the spikes and the reduced system are computed once.
		\item \emph{Guaranteed convergence without tuning.} The iteration converges for every $\lambda>0$ (Theorem~\ref{thm:admm}); there is no relaxation parameter as in PSOR.
		\item \emph{Parallel pointwise steps.} The projection and the multiplier update cost $O(N)$ and are perfectly parallel; the linear solve is distributed by SPIKE.
		\item \emph{One framework.} The same iteration and proof cover the symmetric obstacle problem and the non-symmetric option LCP, and the Black--Scholes and the jump-diffusion models differ only in the matrix and the right-hand side.
		\item \emph{Independence of the partition.} The iterates are identical for every number of processes (Proposition~\ref{prop:spike}); in the experiments the solutions agree to $3.5\times10^{-13}$ or better.
		\item \emph{Scaling in two dimensions.} For $n_x=256$ the time falls from $1343.5$~s to $343.1$~s on four processes (speed-up $3.92$, Section~\ref{sec:perf}).
		\item \emph{Simplicity.} Each iteration is one solve, a $\max$ and an addition; the contact set need not be tracked and nothing is refactorized.
	\end{enumerate}
	
	%======================================================================
	\section{Mathematical theory}\label{sec:theory}
	The properties that the algorithms of Sections~\ref{sec:method} and~\ref{sec:option} rely on are collected here: the existence and uniqueness of the solution, and the structure of the matrices (symmetric positive definite, or diagonally dominant) which makes every diagonal block that SPIKE factorizes invertible. These results are used in the convergence proof of Section~\ref{sec:conv} and in the SPIKE analysis of Section~\ref{sec:spike}.
	\subsection{Well-posedness}
	\begin{proposition}[Unique solution]\label{prop:wp}
		Let $F:\R^N\to\R^N$ be continuous and strongly monotone with modulus $m>0$, $\ip{F(x)-F(y)}{x-y}\ge m\norm{x-y}_2^2$, and let $C\subset\R^N$ be non-empty, closed and convex. Then the variational inequality $u\in C$, $\ip{F(u)}{w-u}\ge0$ for all $w\in C$ has a unique solution. It is equivalent to $-F(u)\in N_C(u)$, the normal cone of $C$ at $u$.
	\end{proposition}
	This is the Lions--Stampacchia / Browder theorem \cite{kinderlehrer}. It covers two cases of interest:
	\begin{itemize}[itemsep=1pt]
		\item \emph{Obstacle problem:} $F(u)=Au-F_0$ with $A$ symmetric positive definite and $C=\{u\ge\varphi\}$. The variational inequality is the optimality condition of $\min_{C}\frac12u^{\top}Au-F_0^{\top}u$, and $m=\lambda_{\min}(A)$. For the five-point matrix $K_h$ on the unit square, $\lambda_{\min}(K_h)=\tfrac{8}{h^2}\sin^2\tfrac{\pi h}{2}\to2\pi^2$.
		\item \emph{Option LCP:} $F(u)=Mu-b$ with $M$ not symmetric but $\ip{x}{Mx}\ge m\norm{x}^2$, $m>0$, and $C=\{u\ge g\}$. The conditions $Mu-b\ge0$, $u\ge g$, $(Mu-b)^{\top}(u-g)=0$ are the same variational inequality. For the log-price scheme the symmetric part of $M$ is $\mathrm{tridiag}(-a,\,2a+c_0,\,-a)$ with $a=\Delta t\,\alpha/h^2$ and $c_0=1+\Delta t(r+\lambda)$, because the convection terms cancel in the symmetric part; hence $\ip{x}{Mx}\ge c_0\norm{x}^2$ and $m=c_0\ge1$. For the Black--Scholes price-variable matrix, the smallest eigenvalue of the symmetric part, evaluated numerically for the matrices of Section~\ref{sec:results}, exceeds $1$.
	\end{itemize}
	
	\subsection{Structure of the matrices}
	\begin{lemma}\label{lem:blocks}
		(a) If $A$ is symmetric positive definite, every principal submatrix of $A+\lambda I$, $\lambda\ge0$, is symmetric positive definite. (b) If $M$ is strictly row diagonally dominant, every principal submatrix of $M+\rho I$, $\rho\ge0$, is strictly row diagonally dominant and hence non-singular. In particular, the diagonal blocks $A_j$ of any partition of such a matrix are invertible.
	\end{lemma}
	\begin{proof}
		(a) For $x\ne0$ supported on the index set, $x^{\top}(A+\lambda I)_{SS}x=\tilde x^{\top}(A+\lambda I)\tilde x>0$ with $\tilde x$ the extension by zeros. (b) Deleting rows and columns does not increase the off-diagonal row sums and leaves the diagonals unchanged. A strictly diagonally dominant matrix is non-singular (Gershgorin). For the log-price matrix, the off-diagonal entries have row sum $|{\rm lower}|+|{\rm upper}|=2\Delta t\,\alpha/h^2<1+\Delta t(2\alpha/h^2+r+\lambda)$ whenever the cell P\'eclet number $|\beta|h/(2\alpha)\le1$, which holds for all grids of Section~\ref{sec:results}.
	\end{proof}
	
	%======================================================================
	\section{Convergence proof}\label{sec:conv}
	We prove convergence of Algorithm~\ref{alg:obst} for the obstacle problem by a Lyapunov argument on the energy of the primal--dual pair $(v^k,d^k)$; the argument holds for every penalty parameter $\lambda>0$ and covers, with the same proof, the time steps of Algorithm~\ref{alg:opt}.
	
	\subsection{Energy functional and splitting}
	Let $K\in\R^{N\times N}$ be symmetric positive definite with $m=\lambda_{\min}(K)>0$ and $L=\lambda_{\max}(K)$ ($K=K_h$ of Section~\ref{sec:method}), let $C=\{v:\ v\ge\varphi\}$, and let
	\begin{equation}\label{eq:energy}
		E(u)=\tfrac12u^{\top}Ku,\qquad u^*=\argmin_{v\in C}E(v).
	\end{equation}
	$E$ is $m$-strongly convex and $C$ is closed and convex, so $u^*$ exists and is unique (Proposition~\ref{prop:wp}). With the indicator $\iota_C$ and the splitting variable $v$, the problem is $\min_{u,v}\ E(u)+\iota_C(v)$ subject to $u=v$, and the scaled ADMM iteration is the alternating minimization of the augmented energy
	\begin{equation}\label{eq:aug}
		\mathcal L_\lambda(u,v;d)=E(u)+\iota_C(v)+\tfrac\lambda2\norm{u-v+d}^2-\tfrac\lambda2\norm d^2,
	\end{equation}
	followed by the multiplier update $d\leftarrow d+u-v$. Written out for the operator $F(u)=Ku$ (more generally for a strongly monotone $F$ with $\ip{F(x)-F(y)}{x-y}\ge m\norm{x-y}^2$):
	\begin{align}
		&u^{k+1}\ \text{the solution of}\ F(u)+\lambda(u-v^k+d^k)=0,\quad\text{i.e.}\ (K+\lambda I)u^{k+1}=\lambda(v^k-d^k), \label{eq:it1}\\
		&v^{k+1}=\Pi_C(u^{k+1}+d^k)=\max(\varphi,\,u^{k+1}+d^k),\qquad d^{k+1}=d^k+u^{k+1}-v^{k+1}. \label{eq:it2}
	\end{align}
	(The penalty parameter is called $\lambda$ here and $\rho$ in Algorithm~\ref{alg:opt}, to avoid confusion with the jump intensity $\lambda$ of the Merton model.) This is Algorithm~\ref{alg:obst} (step~4 is \eqref{eq:it1}, step~5 the projection, step~6 the multiplier update). Equation \eqref{eq:it1} has a unique solution because $F+\lambda I$ is strongly monotone. For the option LCP, $F(u)=Mu-b$ and $\lambda=\rho$. The normal cone of $C$ at $v\in C$ is $N_C(v)=\{y:\ y_i\le0\ \text{for all }i,\ y_i=0\ \text{if }v_i>\varphi_i\}$.
	
	At the solution, $u^*$ satisfies the variational inequality, equivalently $-Ku^*\in N_C(u^*)$ (the discrete complementarity conditions $Ku^*\ge0$, $u^*\ge\varphi$, $(Ku^*)^{\top}(u^*-\varphi)=0$). Put
	\begin{equation}\label{eq:star}
		v^*=u^*,\qquad d^*=-F(u^*)/\lambda,\qquad\text{so}\quad \lambda d^*\in N_C(v^*),\quad F(u^*)=-\lambda d^*.
	\end{equation}
	
	\subsection{The Lyapunov functional and its decrease}
	Define the Lyapunov (energy) functional of the primal--dual pair
	\begin{equation}\label{eq:Vdef}
		V^k=\lambda\norm{d^k-d^*}^2+\lambda\norm{v^k-v^*}^2 .
	\end{equation}
	It measures the distance of $(v^k,d^k)$ from the saddle point $(v^*,d^*)$ of $\mathcal L_\lambda$ in the norm weighted by $\lambda$. Write $r^{k+1}=u^{k+1}-v^{k+1}=d^{k+1}-d^k$ (primal residual) and $s^{k+1}=v^{k+1}-v^k$ (change of the split variable).
	
	\begin{theorem}[Lyapunov decrease and convergence]\label{thm:admm}
		Let $\lambda>0$ and let the iteration \eqref{eq:it1}--\eqref{eq:it2} start from $v^0\in C$ and $d^0$ with $\lambda d^0\in N_C(v^0)$ (for instance $v^0=\varphi$, $d^0=0$, as in Algorithm~\ref{alg:obst}); for arbitrary $v^0,d^0$ the same holds from $k=1$ on. Then for every $k\ge0$ (every $k\ge1$ for arbitrary starting data)
		\begin{equation}\label{eq:lyap}
			V^k-V^{k+1}\ \ge\ 2m\norm{u^{k+1}-u^*}^2+\lambda\norm{r^{k+1}+s^{k+1}}^2\ \ge\ 2m\norm{u^{k+1}-u^*}^2+\lambda\bigl(\norm{r^{k+1}}^2+\norm{s^{k+1}}^2\bigr).
		\end{equation}
		Consequently $V^k$ is non-increasing, $u^k\to u^*$, $v^k\to u^*$, $r^k\to0$, $s^k\to0$ and $d^k\to d^*=-F(u^*)/\lambda$.
	\end{theorem}
	
	\begin{proof}
		\emph{Step 1: relations of the iteration.} From \eqref{eq:it2}, $u^{k+1}+d^k=d^{k+1}+v^{k+1}$. Hence \eqref{eq:it1} gives
		\begin{equation}\label{eq:rel1}
			F(u^{k+1})=-\lambda\,(u^{k+1}-v^k+d^k)=-\lambda\,d^{k+1}-\lambda\,s^{k+1}.
		\end{equation}
		Since $v^{k+1}=\Pi_C(z)$ with $z=u^{k+1}+d^k$, the projection property gives $z-v^{k+1}\in N_C(v^{k+1})$, and $z-v^{k+1}=d^{k+1}$. As $N_C(v)$ is a cone and $\lambda>0$,
		\begin{equation}\label{eq:rel2}
			\lambda\,d^{k+1}\in N_C(v^{k+1}).
		\end{equation}
		
		\emph{Step 2: three inequalities.} Let $e=u^{k+1}-u^*$.
		\begin{enumerate}[label=(\roman*),itemsep=1pt]
			\item Monotonicity of $N_C$, applied to \eqref{eq:rel2} and $\lambda d^*\in N_C(v^*)$: $\ip{d^{k+1}-d^*}{v^{k+1}-v^*}\ge0$.
			\item Strong monotonicity of $F$ with \eqref{eq:rel1} and $F(u^*)=-\lambda d^*$:
			\[
			m\norm e^2\le\ip{F(u^{k+1})-F(u^*)}{e}=-\lambda\ip{d^{k+1}-d^*}{e}-\lambda\ip{s^{k+1}}{e}.
			\]
			\item Monotonicity of $N_C$ applied to \eqref{eq:rel2} at the indices $k+1$ and $k$ (for $k=0$ this uses $\lambda d^0\in N_C(v^0)$): $\ip{r^{k+1}}{s^{k+1}}\ge0$, because $r^{k+1}=d^{k+1}-d^k$ and $s^{k+1}=v^{k+1}-v^k$.
		\end{enumerate}
		
		\emph{Step 3: the energy identity.} Put $a=d^{k+1}-d^*$ and $w=v^{k+1}-v^*$; then $e=r^{k+1}+w$ since $u^*=v^*$. By (i), $\ip ae=\ip a{r^{k+1}}+\ip aw\ge\ip a{r^{k+1}}$, so (ii) gives
		\[
		2m\norm e^2\le-2\lambda\ip a{r^{k+1}}-2\lambda\ip{s^{k+1}}{e}.
		\]
		Since $a=(d^k-d^*)+r^{k+1}$, expanding $\norm a^2$ gives
		\[
		-2\ip a{r^{k+1}}=\norm{d^k-d^*}^2-\norm{d^{k+1}-d^*}^2-\norm{r^{k+1}}^2 .
		\]
		Since $v^k-v^*=w-s^{k+1}$, expanding $\norm{v^k-v^*}^2$ gives $2\ip{s^{k+1}}w=\norm w^2-\norm{v^k-v^*}^2+\norm{s^{k+1}}^2$, and with $e=r^{k+1}+w$,
		\[
		-2\ip{s^{k+1}}{e}=-2\ip{s^{k+1}}{r^{k+1}}+\norm{v^k-v^*}^2-\norm{v^{k+1}-v^*}^2-\norm{s^{k+1}}^2 .
		\]
		Adding the three relations and multiplying by $\lambda$,
		\[
		2m\norm e^2\le V^k-V^{k+1}-\lambda\bigl(\norm{r^{k+1}}^2+\norm{s^{k+1}}^2+2\ip{r^{k+1}}{s^{k+1}}\bigr),
		\]
		which is the first inequality of \eqref{eq:lyap}; the second follows from (iii).
		
		\emph{Step 4: convergence.} By \eqref{eq:lyap}, $V^k\ge0$ is non-increasing, so it converges, and summing \eqref{eq:lyap} gives $\sum_k\bigl(2m\norm{u^{k+1}-u^*}^2+\lambda\norm{r^{k+1}}^2+\lambda\norm{s^{k+1}}^2\bigr)\le V^0<\infty$. Hence $u^k\to u^*$, $r^k\to0$, $s^k\to0$; then $v^{k+1}=u^{k+1}-r^{k+1}\to u^*$, and by \eqref{eq:rel1}, $\lambda d^{k+1}=-F(u^{k+1})-\lambda s^{k+1}\to-F(u^*)=\lambda d^*$.
	\end{proof}
	
	\subsection{Consequences for the obstacle problem}
	\begin{corollary}[Obstacle problem: bounds in energy and in the iterates]\label{cor:obst}
		Let $F(u)=Ku$ and start Algorithm~\ref{alg:obst} from $v^0=\varphi$, $d^0=0$. Then
		\begin{enumerate}[label=(\alph*),itemsep=1pt]
			\item $V^k\le V^0=\dfrac{\norm{Ku^*}^2}{\lambda}+\lambda\norm{\varphi-u^*}^2$ for all $k$, so $\norm{v^k-u^*}^2\le V^0/\lambda$ and $\norm{d^k-d^*}^2\le V^0/\lambda$: the iterates are bounded.
			\item $\displaystyle\min_{0\le j\le k-1}\norm{u^{j+1}-u^*}^2\le\frac{V^0}{2mk}$ for every $k\ge1$.
			\item The energy of the feasible iterate converges to the minimum energy: $\displaystyle\tfrac m2\norm{v^k-u^*}^2\le E(v^k)-E(u^*)\le\norm{Ku^*}\,\norm{v^k-u^*}+\tfrac L2\norm{v^k-u^*}^2\to0$.
		\end{enumerate}
	\end{corollary}
	\begin{proof}
		(a) Theorem~\ref{thm:admm} holds from $k=0$ since $\lambda d^0=0\in N_C(v^0)$; $V^0$ follows from \eqref{eq:Vdef} with $d^*=-Ku^*/\lambda$ and $v^0=\varphi$. (b) Summing \eqref{eq:lyap} over $j=0,\dots,k-1$ gives $2m\sum_{j=0}^{k-1}\norm{u^{j+1}-u^*}^2\le V^0$, and the minimum is at most the mean. (c) $v^k\in C$ and $E(v)-E(u^*)=\ip{Ku^*}{v-u^*}+\tfrac12(v-u^*)^{\top}K(v-u^*)$. The first term is non-negative because $u^*$ solves the variational inequality, and at most $\norm{Ku^*}\norm{v-u^*}$; the second lies between $\tfrac m2\norm{v-u^*}^2$ and $\tfrac L2\norm{v-u^*}^2$. Convergence of $v^k$ is Theorem~\ref{thm:admm}.
	\end{proof}
	
	\begin{remark}[Interpretation]
		(1) The decrease \eqref{eq:lyap} holds with equality up to the non-negative terms dropped in (i) and (iii); it is checked numerically in Section~\ref{sec:results} (the inequality holds at every one of $3000$ iterations for P1). (2) $V^0=\norm{Ku^*}^2/\lambda+\lambda\norm{\varphi-u^*}^2$ is minimal at $\lambda=\norm{Ku^*}/\norm{\varphi-u^*}$, which suggests how the penalty parameter should scale with the data; this is a statement about the upper bound only. (3) The theorem gives convergence of $u^k$ with the non-asymptotic bound (b), of order $O(1/k)$ for the \emph{squared} error of the best iterate so far, and not a linear rate; without strong convexity only an $O(1/k)$ ergodic rate is known in general \cite{hey}. The slow practical convergence in Section~\ref{sec:results} (about $10^3$ iterations in 2D) is consistent with the absence of a proven rate. (4) The stopping rule $\norm{u^{k+1}-u^k}_\infty<\mathrm{tol}$ is controlled by $\norm{r}^2+\norm{s}^2$, which \eqref{eq:lyap} bounds by the decrease of $V$, but a small decrease of $V$ in one step does not bound $V^k$ itself; this is why the 2D errors in Section~\ref{sec:results} are still large when the rule is met. (5) The argument follows \cite[App.~A]{boyd} and \cite{eckstein}; the strong monotonicity of $F$ supplies the extra term $2m\norm{u^{k+1}-u^*}^2$ and hence convergence of $u^k$ itself.
	\end{remark}
	\begin{remark}[Option LCP]
		For $F(u)=Mu-b$ with $\ip{x}{Mx}\ge m\norm x^2$ ($M$ non-symmetric, so $E$ is not an energy but \eqref{eq:lyap} only uses monotonicity), Algorithm~\ref{alg:opt} with $\rho=\lambda$ is the iteration \eqref{eq:it1}--\eqref{eq:it2}, so each time step converges for every $\rho>0$ and the Lyapunov decrease \eqref{eq:lyap} holds with the same $V$. It is the Douglas--Rachford splitting \cite{douglas,lions,eckstein} of the sum of a strongly monotone operator and a normal cone. The theorem is a statement about the exact solution of every inner iteration; the convergence for a fixed number of sweeps is not covered, and Section~\ref{sec:results} shows that it fails in practice.
	\end{remark}
	
	\begin{proposition}[SPIKE exactness and independence of the number of processes]\label{prop:spike}
		Let $B$ be the matrix $A+\lambda I$ of Algorithm~\ref{alg:obst} or $M+\rho I$ of Algorithm~\ref{alg:opt}, partitioned into $p$ consecutive diagonal blocks $B_1,\dots,B_p$ (of any sizes at least $m$, the coupling width, with $m=1$ in 1D and $m=n_x$ in 2D). Then the reduced interface system is non-singular, and the SPIKE solution equals $B^{-1}b$ for every $b$ (in exact arithmetic). Consequently, the iterates $u^k,v^k,d^k$ of Algorithms~\ref{alg:obst} and~\ref{alg:opt} are the same for every $p$ and every partition, and Theorem~\ref{thm:admm} holds for the parallel algorithm.
	\end{proposition}
	\begin{proof}
		By Lemma~\ref{lem:blocks} each diagonal block is invertible, so $D=\mathrm{diag}(B_1,\dots,B_p)$ and $S=D^{-1}B$ are defined; $S$ is non-singular because $B$ is. The system $Sx=g$, $g=D^{-1}b$, consists of the reduced equations for the top and bottom $m$ rows of each block and of the retrieval equations for all rows. If the reduced system had a non-trivial null vector $y$, extending it by the retrieval equations with $g=0$ would give a non-trivial null vector of $S$, a contradiction. So the reduced system has a unique solution and the retrieval equations return the unique solution of $Sx=g$, i.e.\ of $Bx=b$. The map $(v^k,d^k)\mapsto(v^{k+1},d^{k+1})$ of the algorithm depends on the linear solve only through its exact result, so the iterates do not depend on $p$.
	\end{proof}
	
	%======================================================================
	\section{The SPIKE algorithm and its parallel implementation}\label{sec:spike}
	Section~\ref{sec:conv} shows that the ADMM iteration converges whatever solver is used for the linear system, provided that it is exact. We now describe the exact parallel banded solver, SPIKE, that is used for this purpose, its MPI implementation and its cost.
	\subsection{The algorithm}
	We summarize \cite{spike} (see also \cite{polizzi07,sameh78}) for the system $Bx=b$ with $B$ banded. Partition $B$ into $p$ diagonal blocks $B_j$ of order $n_j$ with coupling blocks (above and below the diagonal blocks) of order $m\ll n_j$. Write $B=DS$ with $D=\mathrm{diag}(B_1,\dots,B_p)$; $S$ has identity diagonal blocks and the \emph{spikes}
	\begin{equation}\label{eq:spikes}
		V_j=B_j^{-1}\begin{bmatrix}0\\ X_j\end{bmatrix},\qquad W_j=B_j^{-1}\begin{bmatrix}Y_j\\0\end{bmatrix},
	\end{equation}
	where $X_j$ and $Y_j$ are the coupling blocks to the following and the preceding block. First solve $Dg=b$ (parallel). Then solve $Sx=g$: with the superscripts $(t)$ and $(b)$ for the top and bottom $m$ rows, the interface unknowns satisfy the block-tridiagonal reduced system
	\begin{align}
		x_j^{(t)}+V_j^{(t)}x_{j+1}^{(t)}+W_j^{(t)}x_{j-1}^{(b)}&=g_j^{(t)},\label{eq:red1}\\
		x_j^{(b)}+V_j^{(b)}x_{j+1}^{(t)}+W_j^{(b)}x_{j-1}^{(b)}&=g_j^{(b)},\label{eq:red2}
	\end{align}
	with $V_p=W_1=0$, of dimension $2m(p-1)$ after eliminating the unknowns that do not appear, and the interior values are retrieved independently on each rank,
	\begin{equation}\label{eq:ret}
		x_j=g_j-V_jx_{j+1}^{(t)}-W_jx_{j-1}^{(b)}.
	\end{equation}
	For strictly diagonally dominant matrices a truncated version, in which the spikes are cut off away from the interface so that the reduced system becomes block diagonal, is available \cite{spike}. Our implementation builds the full spikes and solves the full reduced system.
	
	\subsection{MPI structure}
	The code is written in Fortran~90 with MPI \cite{mpi}. Rank $j$ owns the $n/p$ consecutive rows of $B_j$ and the matching parts of the vectors. The setup of Algorithms~\ref{alg:obst} and~\ref{alg:opt} is done once: (i) factor $B_j$; (ii) compute the spikes (the interface parts $V_j^{(t,b)},W_j^{(t,b)}$ and the interior parts); (iii) gather the interface blocks on rank~0, assemble the reduced matrix and factor it. Each ADMM iteration then performs:
	\begin{enumerate}[itemsep=1pt]
		\item rank~0 forms the right-hand side from $v^k,d^k$ and broadcasts it (projection and multiplier update are $O(n)$ and are done on rank~0);
		\item each rank solves $B_jg_j=b_j$ and sends the $2m$ interface entries of $g_j$ to rank~0;
		\item rank~0 solves the reduced system and broadcasts the $2m(p-1)$ interface values;
		\item each rank recovers its interior values by \eqref{eq:ret} and sends them back to rank~0.
	\end{enumerate}
	In 2D a row-slab partition needs every slab to contain at least two grid rows, so $p\le n_x/2$.
	
	\subsection{Cost model}
	With a banded Cholesky or LU factorization of the diagonal blocks (half-bandwidth $b=1$ in 1D, $b=n_x$ in 2D), the setup costs to leading order $O\bigl((n/p)\,b^2\bigr)$ for the factorization and $O\bigl((n/p)\,b\bigr)$ per spike column, i.e.\ $O\bigl((n/p)\,b\,m\bigr)=O((n/p)b^2)$ for the $2m$ spike columns in 2D; each local solve costs $O((n/p)\,b)$; the reduced system has dimension $2m(p-1)$ and costs $O(p\,m^3)$ if factorized by a block-tridiagonal recursion. These estimates are ours and are not measured; they explain the qualitative behaviour of Table~\ref{tab:t2d} (the part $O(n_x^4/p)$ scales, the reduced solve does not) but are not fitted to it.
	
	%======================================================================
	\section{Parallel performance}\label{sec:perf}
	This section reports the cost of the parallel implementation of Section~\ref{sec:spike}; the accuracy of the computed solutions is the subject of Section~\ref{sec:results}. The timings use the obstacle problems that are defined in Section~\ref{sec:problems}.
	\subsection{Measurements}
	Table~\ref{tab:t1d} gives timings of the MPI SPIKE-ADMM implementation for the 1D obstacle problem and Table~\ref{tab:t2d} for the 2D problem. Every process reports its CPU time in seconds (measured inside the program); we report the maximum over processes. % TODO: add machine, compiler and MPI library here.
	The 1D obstacle runs use $50$ ADMM iterations. In 2D, $n_x$ is the number of grid points per direction ($n_x^2$ unknowns). A dash means that the run was not made (for 2D, $P\le n_x/2$).
	\begin{table}[t]
		\centering\small
		\caption{1D obstacle problem: maximum over processes of the reported time (s), $n$ grid points, $P$ processes.}\label{tab:t1d}
		\begin{tabular}{rrrrrr}
			\toprule
			$n$ & $P=1$ & $P=2$ & $P=4$ & $P=8$ & $P=16$\\
			\midrule
			16 & 0.0043 & 0.0012 & 0.0012 & 0.0022 & -- \\
			32 & 0.0175 & 0.0026 & 0.0023 & 0.0024 & 0.0081 \\
			64 & 0.0698 & 0.0114 & 0.0028 & 0.0030 & 0.0082 \\
			128 & 0.4922 & 0.0765 & 0.0174 & 0.0050 & 0.0125 \\
			256 & 5.355 & 0.5497 & 0.0901 & 0.0245 & 0.0151 \\
			512 & 45.38 & 3.376 & 0.3698 & 0.0803 & 0.0202 \\
			1024 & 530.3 & 94.02 & 5.963 & 1.403 & 0.2071 \\
			2048 & -- & 574.1 & 57.19 & 13.44 & 1.922 \\
			4096 & -- & -- & 690.2 & 145.9 & 16.11 \\
			8192 & -- & -- & -- & 1140.5 & 200.0 \\
			\bottomrule
		\end{tabular}
	\end{table}
	\begin{table}[t]
		\centering\small
		\caption{2D obstacle problem: maximum over processes of the reported time (s), $n_x\times n_x$ grid.}\label{tab:t2d}
		\begin{tabular}{rrrrrr}
			\toprule
			$n_x$ & $P=1$ & $P=2$ & $P=4$ & $P=8$ & $P=16$\\
			\midrule
			16 & 0.0176 & 0.0134 & 0.0078 & 0.0057 & -- \\
			32 & 0.1719 & 0.0957 & 0.0498 & 0.0343 & 0.0190 \\
			64 & 2.592 & 1.292 & 0.6556 & 0.5354 & 0.4424 \\
			256 & 1343.5 & 682.5 & 343.1 & 248.0 & 245.8 \\
			\bottomrule
		\end{tabular}
	\end{table}
	\begin{figure}[t]
		\centering
		\includegraphics[width=\textwidth]{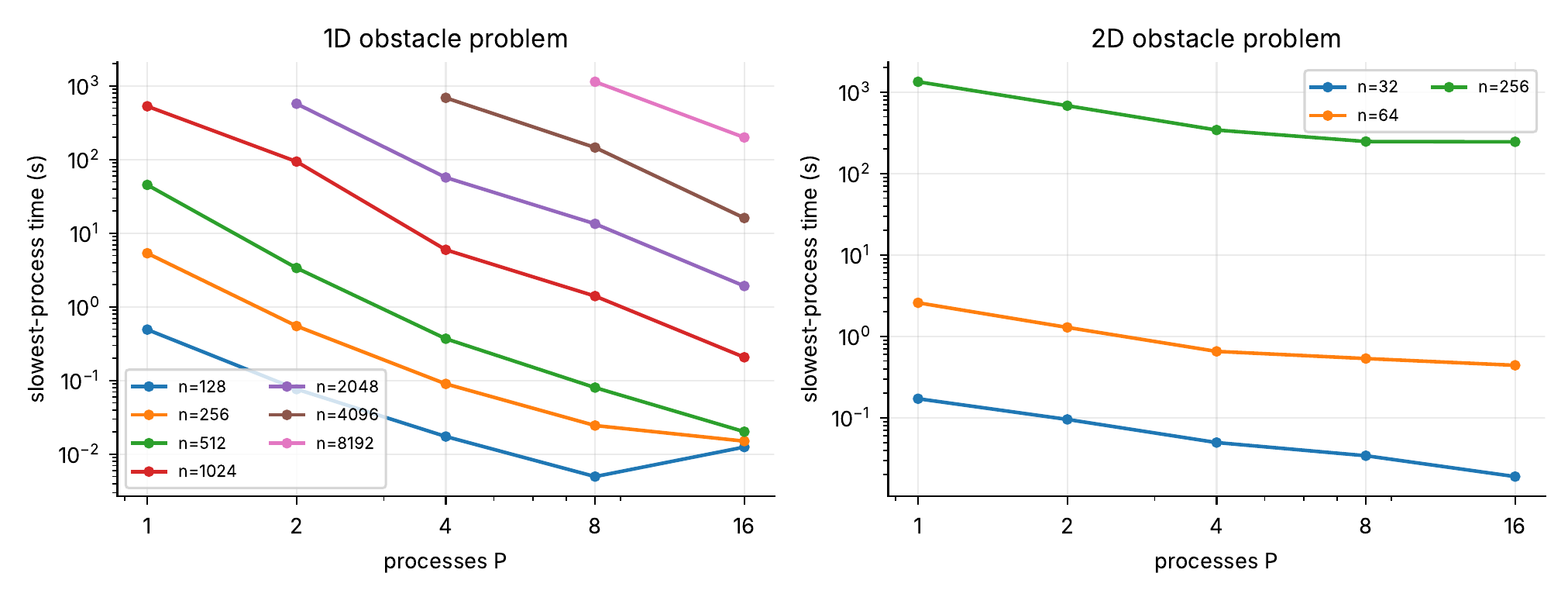}
		\caption{Maximum per-process time against the number of processes, from Tables~\ref{tab:t1d} and~\ref{tab:t2d} (log--log axes).}
		\label{fig:fig1timing2dspeedup}
	\end{figure}
	
	\subsection{Observations and limits of the data}
	\textbf{1D.} With $P=1$, doubling $n$ multiplies the time by $7.1$, $10.9$ and $8.5$ from $n=64$ to $512$ (and by $11.7$ for $512\to1024$), i.e.\ roughly cubic growth, which is what dense block work $O((n/p)^3)$ would give and not what a banded algorithm would; for the same reason the times fall faster than $1/P$ ($n=1024$: $530.3$~s on one, $94.0$~s on two and $5.96$~s on four processes). In the code used for the 1D timings the diagonal blocks are treated as dense matrices (a dense inverse of every block), so the work per block is $O((n/p)^3)$ instead of $O(n/p)$. A speed-up computed against the $P=1$ run therefore exceeds $P$ and measures the change in block size, not parallel efficiency; we do not report such speed-ups as a result. For small $n$ the runs with more processes are slower ($n=16$: $P=8$ against $P=4$; $n=128$: $P=16$ against $P=8$) because communication and the reduced solve dominate.
	
	\textbf{2D.} The $P=1$ times for $n_x=16,32,64$ grow by factors $9.8$ and $15.1$ per doubling, between $n_x^3$ and $n_x^4$, which matches a banded factorization (bandwidth squared times the number of unknowns). For $n_x=256$ the times are $1343.5$, $682.5$, $343.1$, $248.0$ and $245.8$~s for $P=1,2,4,8,16$, i.e.\ speed-ups $1.97$, $3.92$, $5.42$ and $5.47$, with parallel efficiencies $98\%$, $98\%$, $68\%$ and $34\%$: almost linear up to $P=4$ and no gain from $8$ to $16$ processes. We have not measured the cause of the saturation; candidates are the serial reduced solve on rank~0 (its dimension $2n_x(P-1)$ grows with $P$), data gathered on rank~0, and memory bandwidth. For $n_x=64$ the speed-ups are $2.0$, $4.0$, $4.8$ and $5.9$ for $P=2,4,8,16$, i.e.\ efficiencies of $100\%$, $99\%$, $61\%$ and $37\%$.
	
	\textbf{Baseline missing.} For a tridiagonal system SPIKE cannot be expected to beat a serial Thomas solve; for the 2D problem a serial banded-LU ADMM run on the same machine is the proper baseline. Parallel claims in this paper are therefore limited to the scaling of the SPIKE solver itself.
	
	%======================================================================
	\section{Numerical problems}\label{sec:problems}
	The numerical problems are seven: three obstacle problems from \cite{tran} (Section~6, Eqs.~(6.1), (6.2), (6.4)) and four option-pricing problems. They are stated below together with the data, the discretization and the reference used (Table~\ref{tab:probs}). The same ADMM--SPIKE solver (Algorithms~\ref{alg:obst} and~\ref{alg:opt}) is used for all of them; only the matrix and the right-hand side change.
	
	\subsection{Obstacle problems (P1--P3)}
	\paragraph{Problem 1: one-dimensional obstacle $\varphi_1$ \cite[Eq.~(6.1)]{tran}.} Find $u=\argmin_{v\ge\varphi_1}\tfrac12\int_0^1|v'|^2dx$, $v(0)=v(1)=0$, with
	\begin{equation}\label{eq:P1}
		\varphi_1(x)=\begin{cases}100x^2,& 0\le x\le0.25,\\ 100x(1-x)-12.5,& 0.25\le x\le0.5,\\ \varphi_1(1-x),& 0.5\le x\le1.0.\end{cases}
	\end{equation}
	The exact solution is $u_1$ of \eqref{eq:exact}. Data: $256$ interior grid points ($h=1/257$), penalty parameter $\lambda=45$, $50$ iterations.
	
	\paragraph{Problem 2: one-dimensional obstacle $\varphi_2$ \cite[Eq.~(6.2)]{tran}.} The same problem with the trigonometric obstacle
	\begin{equation}\label{eq:P2}
		\varphi_2(x)=\begin{cases}10\sin(2\pi x),& 0\le x\le0.25,\\ 5\cos\bigl(\pi(4x-1)\bigr)+5,& 0.25\le x\le0.5,\\ \varphi_2(1-x),& 0.5\le x\le1.0.\end{cases}
	\end{equation}
	The exact solution is $u_2$ of \eqref{eq:exact}. Data: $256$ interior grid points, $\lambda=250$ (the value used in \cite{tran} for this obstacle).
	
	\paragraph{Problem 3: two-dimensional obstacle, non-smooth \cite[Eq.~(6.4)]{tran}.} Find $u=\argmin_{v\ge\varphi_3}\tfrac12\int_\Omega|\nabla v|^2dx$ on $\Omega=(0,1)^2$, $v=0$ on $\partial\Omega$, with the non-smooth obstacle
	\begin{equation}\label{eq:P3}
		\varphi_3(x,y)=\begin{cases}5.0,&|x-0.6|+|y-0.6|<0.04,\\ 4.5,&(x-0.6)^2+(y-0.25)^2<0.001,\\ 4.5,&y=0.57\ \text{and}\ 0.075<x<0.13,\\ 0,&\text{otherwise.}\end{cases}
	\end{equation}
	No closed-form solution is known. Data: $n_x\times n_x$ interior grid, $h=1/(n_x+1)$, $n_x=64,128,256$ for the accuracy runs and $n_x=16,32,64,256$ for the timings; $\lambda=1.3\times10^4$ and stopping rule $\|u^{k+1}-u^k\|_\infty<5\times10^{-4}$ as in \cite{tran}.
	
	\subsection{Option-pricing problems (P4--P7)}
	\paragraph{Problem 4: Black--Scholes call in the transformed form.} With the shift $\varphi=u-\psi$ the option problem on $(-\bar x,\bar x)$ takes the form
	\begin{align}
		\frac{\partial\varphi}{\partial t}-\frac12\sigma^2(x,t)\frac{\partial^2\varphi}{\partial x^2}-\beta(x,t)\frac{\partial\varphi}{\partial x}+r(t)\,\varphi&=f(x,t),\quad x\in(-\bar x,\bar x),\ t\in(0,T],\label{eq:P7a}\\
		\varphi(-\bar x,t)=\varphi(\bar x,t)&=0,\quad t\in(0,T],\label{eq:P7b}\\
		\varphi(x,0)&=\varphi_0(x),\quad x\in(-\bar x,\bar x),\label{eq:P7c}
	\end{align}
	where $\sigma(x,t)$, $\beta(x,t)$ and $r(t)$ are the coefficients of the transformed problem, $f$ and $\varphi_0$ are the source term and the initial value that follow from the shift, $\bar x$ is the truncation point of the log-price, $K=E$ is the strike, and $\psi(x,t)=Ke^{-\int_0^tr(\tau)\,d\tau}-e^{x}$ in $[-\infty,-\bar x]\times[0,T]$, $\psi(x,t)=0$ in $[x_0,\infty]\times[0,T]$, and $\psi$ is smooth. The boundary data are thereby removed from the discrete system. Data: Black--Scholes call, $K=100$, $r=0.1$, $\sigma=0.2$, $T=1$, $S\in[0,200]$, $V(0)=0$, $V(200)=200-Ke^{-rT}$; $512$ intervals, $1000$ time steps, $\rho=1$ and tolerance $10^{-10}$. Without dividends early exercise of a call is never optimal \cite{merton73}, so the American and the European value coincide and the closed-form Black--Scholes value is the reference.
	
	\paragraph{Problem 5: American put under Black--Scholes.} With $\tau$ the time to maturity and $g$ the payoff,
	\begin{equation}\label{eq:P4}
		u(\tau,x):=\begin{cases}u_\tau(\tau,x)-\mathcal Lu(\tau,x)\ge0,\\ u(\tau,x)\ge g(x),\\ \bigl(u_\tau(\tau,x)-\mathcal Lu(\tau,x)\bigr)\bigl(u(\tau,x)-g(x)\bigr)=0,\end{cases}\qquad (\tau,x)\in(0,T)\times(-\infty,\infty).
	\end{equation}
	For the Black--Scholes model $\mathcal Lu=\frac12\sigma^2S^2u_{SS}+rSu_S-ru$. The American put has $g(S)=\max(K-S,0)$ with $K=100$, $r=0.05$, $\sigma=0.2$, $T=1$ and $S\in[0,200]$; Dirichlet values at the ends are the payoff values $V(0)=K$ and $V(200)=0$. Discretization: uniform grid in $S$ with $511$ interior nodes ($\Delta S\approx0.39$), central differences, implicit Euler, $1000$ time steps; each step is the discrete LCP \eqref{eq:dlcp} with a tridiagonal non-symmetric matrix, solved by Algorithm~\ref{alg:opt} with $\rho=1$ and tolerance $10^{-10}$. Reference: binomial tree \cite{crr}, $5000$ steps.
	
	\paragraph{Problem 6: American put under Merton's jump-diffusion model.} In the log-price $x$ the operator is
	\begin{equation}\label{eq:P6a}
		\mathcal Lu(x,t)=\mathcal Du(x,t)+\lambda_J\,\mathcal Iu(x,t),
	\end{equation}
	where $\mathcal D$ and $\mathcal I$ are the differential and integral operators
	\begin{align}
		\mathcal Du(x,t)&=\tfrac12\sigma^2\frac{\partial^2u}{\partial x^2}(x,t)+\Bigl(r-\frac{\sigma^2}{2}-\lambda_J\kappa\Bigr)\frac{\partial u}{\partial x}(x,t)-(r+\lambda_J)\,u(x,t),\label{eq:P6b}\\
		\mathcal Iu(x,t)&=\int_{-\infty}^{\infty}u(y,t)\,\gamma(y-x)\,dy.\label{eq:P6c}
	\end{align}
	To approximate $\mathcal I$ numerically, the integral is separated into its parts on $\Omega$ and on $\mathbb R\setminus\Omega$,
	\begin{equation}\label{eq:P6d}
		\mathcal Iu(x,t)=\int_\Omega u(y,t)\,\gamma(y-x)\,dy+\mathcal R^*(x,t),\qquad \mathcal R^*(x,t)=\int_{\mathbb R\setminus\Omega}u(y,t)\,\gamma(y-x)\,dy,
	\end{equation}
	where, for a European put, the asymptotic form gives
	\begin{equation}\label{eq:P6e}
		\mathcal R^*(x,t)=\int_{\mathbb R\setminus\Omega}\bigl(Ee^{-rt}-Ee^{y}\bigr)\,\gamma(y-x)\,dy .
	\end{equation}
	The drift contains the compensator $\kappa=\int(e^{y}-1)\gamma(y)\,dy=e^{\mu_J+\sigma_J^2/2}-1$ for Merton's model. \emph{Data (the parameters of \cite{maurya}).} $\sigma=0.15$, $r=0.05$, jump intensity $\lambda_J=0.10$, $\mu_J=-0.90$, $\sigma_J=0.45$, $E=100$, $T=0.25$, so $\kappa\approx-0.550$; the quantity computed is the American put with the problem \eqref{eq:P4} and $\mathcal L$ of \eqref{eq:P6a}. \emph{Discretization.} Truncated domain $(x_{\min},x_{\max})=(-1.5,1.5)$, $M=800$ intervals, $N=400$ time steps; the differential part is implicit (tridiagonal matrix, SPIKE), the integral term is evaluated explicitly by the trapezoidal rule at the previous time level, and the tail is added analytically (Problem~7); Algorithm~\ref{alg:opt} with $\rho=1$, tolerance $10^{-10}$. \emph{Reference.} The American put prices at $S=90,100,110$ of \cite{maurya}.
	
	\paragraph{Problem 7: European options and the tails of the integral operator.} Evaluated for the two jump densities $\gamma$, the tail $\mathcal R^*$ of \eqref{eq:P6d}--\eqref{eq:P6e} is
	\begin{equation}\label{eq:P6f}
		\small\mathcal R^*(x,t)=\begin{cases}
			E(1-p)\,e^{\zeta_2(x_{\min}-x)-rt}-E(1-p)\dfrac{\zeta_2}{\zeta_2+1}\,e^{-\zeta_2x+(\zeta_2+1)x_{\min}}&\text{Kou's model,}\\[2mm]
			Ee^{-rt}\,\Phi\!\Bigl(\dfrac{x_{\min}-x-\mu_J}{\sigma_J}\Bigr)-Ee^{x+\mu_J+\frac{\sigma_J^2}{2}}\,\Phi\!\Bigl(\dfrac{x_{\min}-x-\mu_J-\sigma_J^2}{\sigma_J}\Bigr)&\text{Merton's model,}
		\end{cases}
	\end{equation}
	where $\Phi(z)=\frac1{\sqrt{2\pi}}\int_{-\infty}^{z}e^{-s^2/2}\,ds$ is the distribution function of the standard Gaussian. (In the last argument of the Merton case we use $\sigma_J^2$, which is what the derivation gives and what \cite{maurya} has.) In Kou's model $\gamma(y)=p\,\zeta_1e^{-\zeta_1y}\mathbf 1_{y>0}+(1-p)\,\zeta_2e^{\zeta_2y}\mathbf 1_{y<0}$ with $p$ the probability of an upward jump and $\zeta_1,\zeta_2>0$ the decay rates of the up and down jumps. Only Merton's model is treated numerically; Kou's tail is stated for completeness. \emph{Data.} The parameters of Problem~6; the options priced are the European put (tail \eqref{eq:P6f}, Merton case) and the European call (the analogue \eqref{eq:tailC} for the right tail), both with $E=100$, $T=0.25$. \emph{Reference.} The prices at $S=90,100,110$ of \cite{maurya} and Merton's Poisson-mixture series \cite{merton}.
	
	\begin{table}[t]
		\centering\footnotesize
		\caption{Summary of the seven numerical problems. $\lambda$ in P1--P3 is the ADMM penalty parameter of Algorithm~\ref{alg:obst}; $\lambda_J=0.10$ in P6, P7 is the jump intensity.}\label{tab:probs}
		\setlength{\tabcolsep}{3pt}
		\begin{tabular}{@{}p{0.9cm}p{3.1cm}p{3.1cm}p{3.0cm}p{2.6cm}@{}}
			\toprule
			& Problem & Data & Discretization & Reference\\
			\midrule
			P1 & 1D obstacle $\varphi_1$ \eqref{eq:P1} & $\lambda=45$ & $256$ interior points, $h=1/257$ & exact $u_1$ \eqref{eq:exact}\\
			P2 & 1D obstacle $\varphi_2$ \eqref{eq:P2} & $\lambda=250$ & $256$ interior points, $h=1/257$ & exact $u_2$ \eqref{eq:exact}\\
			P3 & 2D obstacle $\varphi_3$ \eqref{eq:P3} & $\lambda=1.3\times10^4$, tol $5\times10^{-4}$ & $n_x=64,128,256$; five-point stencil & discrete solution of an active-set solver with SPIKE solves\\
			P4 & American option, Example 3 \eqref{eq:P7a}--\eqref{eq:P7c} (Black--Scholes call) & $K=100$, $r=0.1$, $\sigma=0.2$, $T=1$, $S\in[0,200]$ & $512$ intervals, $1000$ steps & closed-form Black--Scholes\\
			P5 & American put, Black--Scholes \eqref{eq:P4} & $K=100$, $r=0.05$, $\sigma=0.2$, $T=1$, $S\in[0,200]$ & $511$ nodes, $1000$ steps & binomial tree, $5000$ steps\\
			P6 & Merton jump-diffusion PIDE \eqref{eq:P6a}--\eqref{eq:P6e}, American put & $\sigma=0.15$, $r=0.05$, $\lambda_J=0.10$, $\mu_J=-0.90$, $\sigma_J=0.45$, $E=100$, $T=0.25$ & $x\in(-1.5,1.5)$, $M=800$, $N=400$ & American put of \cite{maurya}\\
			P7 & Tails \eqref{eq:P6f}: European put and call under Merton & as P6 & as P6 & prices of \cite{maurya}; Merton series \cite{merton}\\
			\bottomrule
		\end{tabular}
	\end{table}
	
	\subsection{Timing experiments}
	Separately from P1--P7, the MPI programs were timed for the 1D obstacle problem ($n=16,\dots,8192$ points, $50$ ADMM iterations) and the 2D obstacle problem ($n_x=16,32,64,256$ points per direction) on $P=1,2,4,8,16$ processes, with $P\le n_x/2$ in 2D; the maximum over processes of the reported CPU time is given in Tables~\ref{tab:t1d} and~\ref{tab:t2d}.
	
	%======================================================================
	\section{Numerical results}\label{sec:results}
	The test problems of Section~\ref{sec:problems} are solved by the implementation of Section~\ref{sec:spike}. The one- and two-dimensional obstacle problems P1--P3 verify the method for the symmetric case, and the option problems P4--P7 verify it for the non-symmetric case, in the order of increasing difficulty: Black--Scholes call and put, and then Merton's jump-diffusion model.
	\subsection{Experimental setup}
	The ADMM--SPIKE solver of Sections~\ref{sec:method}--\ref{sec:spike} is implemented in Fortran~90 with MPI in double precision: the diagonal blocks, the spikes and the reduced system are factorized once, and every iteration consists of the steps of Section~\ref{sec:spike}. The accuracy is measured against exact solutions (P1, P2), closed-form and binomial prices (P4, P5), published prices and Merton's series (P6, P7), and, for P3, against an active-set solver. Solutions computed with different numbers of processes $p$ are compared with each other to verify Proposition~\ref{prop:spike}.
	
	\subsection{P1, P2: one-dimensional obstacle problems}
\begin{figure}[t]
		\centering
		\includegraphics[width=0.7\linewidth]{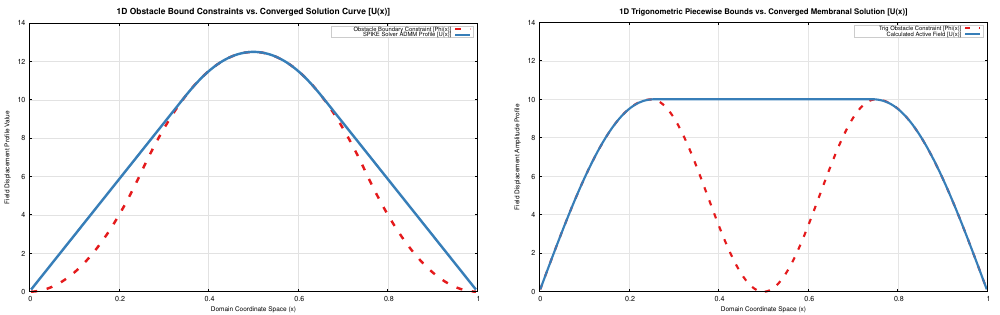}
		\caption{P1 (left, $\varphi_1$) and P2 (right, $\varphi_2$). Dashed: obstacle; solid: ADMM--SPIKE solution.}
		\label{fig:fig21dobstaclep1p2}
	\end{figure}

\begin{figure}[t]
		\centering
		\includegraphics[width=0.7\linewidth]{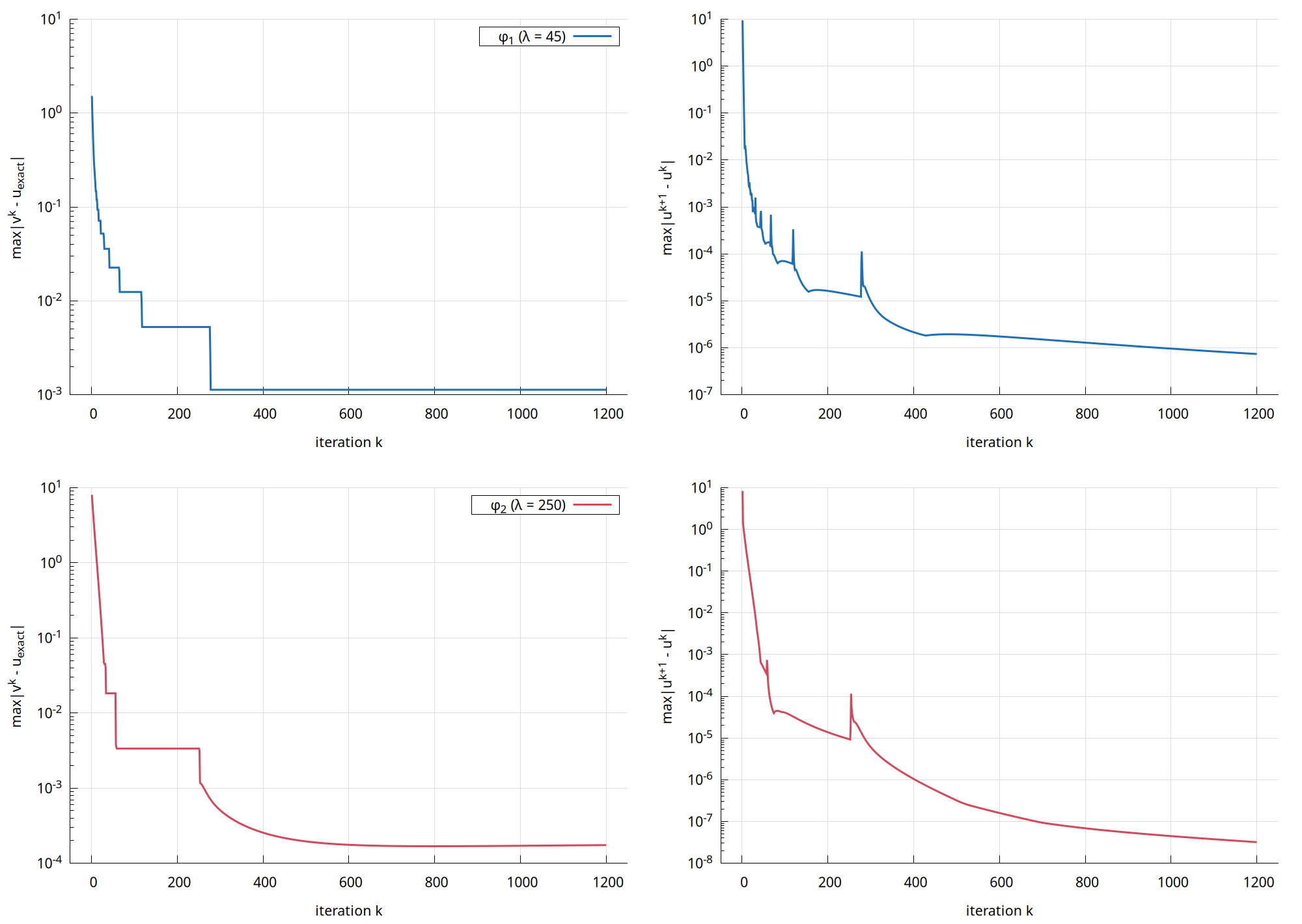}
		\caption{Convergence of the ADMM--SPIKE iteration for P1 (top, $\varphi_1$, $\lambda=45$) and P2 (bottom, $\varphi_2$, $\lambda=250$): maximum error $\max|v^k-u_{\mathrm{exact}}|$ of the feasible iterate (left) and change per iteration $\max|u^{k+1}-u^k|$ (right) against the iteration number $k$ ($256$ interior nodes).}
		\label{fig:fig3convergencephi1phi2}
	\end{figure}

	Figure~\ref{fig:fig21dobstaclep1p2} shows the computed solutions. The solution (solid) coincides with the obstacle (dashed) on the contact set, which is the neighbourhood of the maximum for $\varphi_1$ and the two outer arcs for $\varphi_2$; between them it is a straight line (for $\varphi_1$, of slope $100-50\sqrt2\approx29.3$) or the plateau $10$ (for $\varphi_2$). Both are the exact solutions \eqref{eq:exact}.
	Table~\ref{tab:1d} gives the maximum error of the feasible iterate $v^k$ against the exact solution after $k$ iterations ($256$ interior nodes, $h=1/257$); the error of $u^k$ is of the same order ($5.4\times10^{-2}$, $8.6\times10^{-3}$, $4.3\times10^{-3}$, $2.1\times10^{-3}$, $8.9\times10^{-4}$ for $\varphi_1$ at $k=10,50,100,200,400$). The error falls slowly with $k$: after $50$ iterations it is $2.3\times10^{-2}$ for $\varphi_1$ and $1.8\times10^{-2}$ for $\varphi_2$. Run to $2\times10^4$ iterations the error for $\varphi_1$ reaches $2.8\times10^{-5}$; for $\varphi_2$ the change per iteration falls below $10^{-12}$ after $10084$ iterations with an error of $1.9\times10^{-4}$, which is the discretization error of the grid (the kink at $x=0.25$ is not a grid node). The maximum error can stay constant over a range of iterations ($3.4\times10^{-3}$ for $\varphi_2$ at $k=100$ and $k=200$); this is consistent with the maximum being attained at one node while the iterate in the contact region changes slowly (Figure~\ref{fig:fig3convergencephi1phi2}). The solutions computed with $p=2,4,8$ processes agree with the $p=1$ solution after $50$ iterations to $6.0\times10^{-14}$ ($\varphi_1$) and $1.8\times10^{-13}$ ($\varphi_2$), as Proposition~\ref{prop:spike} predicts. The Lyapunov decrease \eqref{eq:lyap} was checked for P1 ($\lambda=45$, $v^0=\varphi_1$, $d^0=0$, $u^*$ the exact discrete solution computed by an active-set method): it holds at each of $3000$ iterations, $V^k$ is monotonically decreasing ($V^0=8.33\times10^4$, $V^{10}=3.68\times10^3$, $V^{100}=8.28\times10^2$, $V^{1000}=48.9$, $V^{3000}=5.15$), and the closed form $V^0=\norm{Ku^*}^2/\lambda+\lambda\norm{\varphi-u^*}^2$ of Corollary~\ref{cor:obst} agrees with the computed value.
	\begin{table}[t]
		\centering\small
		\caption{P1, P2: maximum error of the feasible iterate $v^k$ against the exact solution \eqref{eq:exact} after $k$ ADMM iterations ($256$ interior nodes).}\label{tab:1d}
		\begin{tabular}{lrrrrr}
			\toprule
			Obstacle & $k=10$ & $k=50$ & $k=100$ & $k=200$ & $k=400$\\
			\midrule
			$\varphi_1$ ($\lambda=45$) & $1.5\times10^{-1}$ & $2.3\times10^{-2}$ & $1.2\times10^{-2}$ & $5.3\times10^{-3}$ & $1.1\times10^{-3}$\\
			$\varphi_2$ ($\lambda=250$) & $1.6$ & $1.8\times10^{-2}$ & $3.4\times10^{-3}$ & $3.4\times10^{-3}$ & $2.5\times10^{-4}$\\
			\bottomrule
		\end{tabular}
	\end{table}
	
	\subsection{P3: two-dimensional obstacle problem}
\begin{figure}[t]
		\centering
		\includegraphics[width=\textwidth]{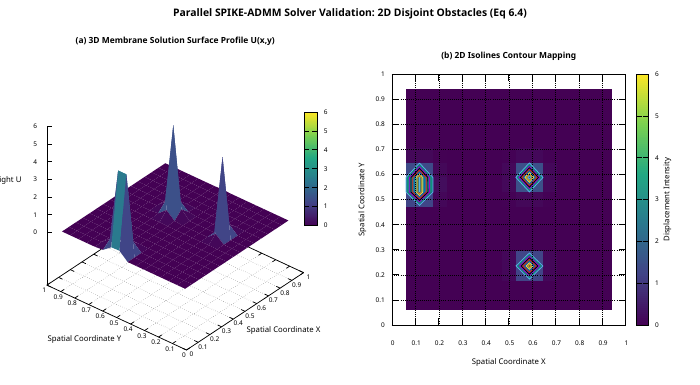}
		\caption{P3: (a) surface of the membrane solution; (b) iso-lines.}
		\label{fig:fig4p3obstacle2d}
	\end{figure}

	Figure~\ref{fig:fig4p3obstacle2d} shows the computed solution: three peaks at the diamond, the disc and the segment of $\varphi_3$, and a flat solution elsewhere.
	For $\lambda=1.3\times10^4$ the iterates computed with $p=1$ and $p=4$ slabs agree to $9.1\times10^{-14}$ ($n_x=64$, $300$ iterations), $1.6\times10^{-13}$ ($n_x=128$, $300$ iterations) and $2.8\times10^{-13}$ ($n_x=256$, $2000$ iterations). The change per iteration falls below $5\times10^{-4}$ after $1038$, $1032$ and $1030$ iterations for $n_x=64,128,256$ (i.e.\ about $10^3$ on every grid). Compared with the solution of an independent active-set solver (an active-set solver with SPIKE linear solves, cf.\ \cite{singh}, run to the KKT tolerance $10^{-6}$), the error of the ADMM iterate at that point is however still large: $0.18$ ($n_x=64$), $0.18$ ($n_x=128$) and $0.18$ ($n_x=256$). After $5000$ iterations the error is $3.0\times10^{-6}$ ($n_x=64$) and $2.5\times10^{-6}$ ($n_x=128$), and after $2000$ iterations $1.1\times10^{-2}$ ($n_x=256$). The stopping rule $\norm{u^{k+1}-u^k}_\infty<5\times10^{-4}$ of \cite{tran} is therefore not a reliable measure of the error for this obstacle in the hard-constraint form; a residual-based or KKT-based test (as in \cite{tran}'s error plots,) is preferable. The active-set hybrid needed $2$, $5$ and $9$ outer iterations for the same three grids.
	
	\subsection{P4: Black--Scholes call}
	For $K=100$, $r=0.1$, $\sigma=0.2$, $T=1$ the ADMM--SPIKE values of the call at $S=80,90,100,110,120$ are $2.7909$, $6.9487$, $13.2682$, $21.2476$, $30.2576$, against the closed-form values $2.7899$, $6.9490$, $13.2697$, $21.2488$, $30.2585$ (Figure~\ref{fig:fig5blackscholesputcall}, right); the error is at most $1.5\times10^{-3}$, which is the size of the first-order time-stepping error.
	
	\subsection{P5: Black--Scholes American put}\label{sec:p5}
\begin{figure}[t]
		\centering
		\includegraphics[width=\textwidth]{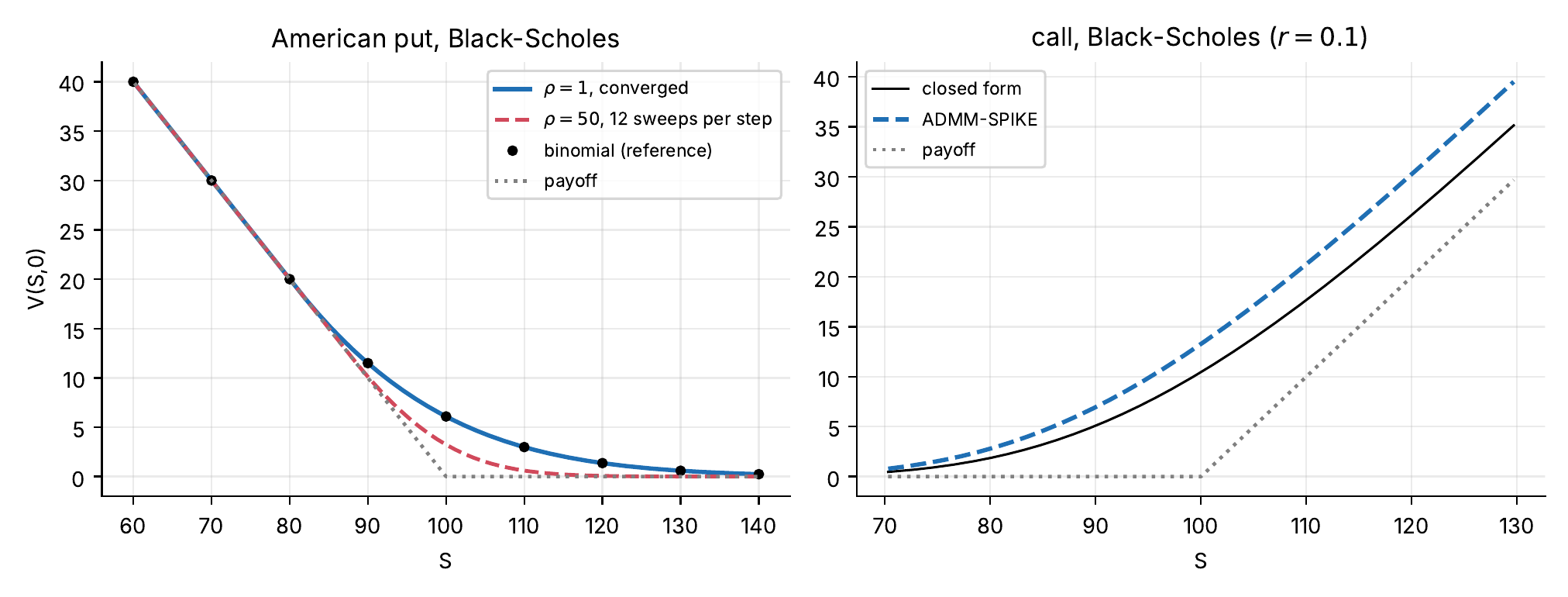}
		\caption{P5 (left) and P4 (right). Left: American put, ADMM-SPIKE with $\rho=1$ converged, with $\rho=50$ and $12$ sweeps per step, and binomial reference values. Right: Black--Scholes call, closed form and ADMM-SPIKE.}
		\label{fig:fig5blackscholesputcall}
	\end{figure}

	The independent references are a binomial tree \cite{crr} with $5000$ steps, $V(90)=11.4928$, $V(100)=6.0902$, $V(110)=2.9866$, and, for the European put, the closed-form value $5.5735$ at $S=100$. With $511$ interior nodes, $1000$ time steps and $p=4$, the ADMM with $\rho=1$, run to $10^{-10}$ in every step (on average $55$ sweeps per step), gives $V=20.0000$, $11.4919$, $6.0885$, $2.9853$, $1.3665$ at $S=80,90,100,110,120$, within $1.8\times10^{-3}$ of the binomial values (Figure~\ref{fig:fig5blackscholesputcall}, left). The $p=4$ and $p=1$ solutions agree to $6.1\times10^{-14}$. If the inner iteration is stopped after a fixed number of sweeps, the result depends strongly on $\rho$ (Table~\ref{tab:rho}): with $\rho=50$ and $12$ sweeps per step the values are $20.0000$, $10.1119$, $3.2344$, $0.6046$, $0.0662$, i.e.\ $V(100)=3.23$, an error of $47\,\%$, which is the under-converged ADMM of Section~\ref{sec:option}. The error depends on the augmentation parameter relative to the diagonal of $M$ ($1.0$ to $11.4$) and on the sweep budget: with $\rho=1$ and the same $12$ sweeps it is $0.04\,\%$, and every value of $\rho$ gives $0.03\,\%$ when the inner iteration is run to convergence, at a cost of $42$ to $691$ sweeps per step.
	\begin{table}[t]
		\centering\small
		\caption{P5: relative error (\%) of the American put at $S=100$ against the binomial value $6.0902$ ($5000$ steps), as a function of the augmentation parameter $\rho$ and the number of ADMM sweeps per time step ($511$ nodes, $1000$ steps, $p=4$). The diagonal of $M$ ranges from $1.0$ to $11.4$.}\label{tab:rho}
		\begin{tabular}{rrrrr}
			\toprule
			$\rho$ & 6 sweeps & 12 sweeps & converged ($10^{-10}$) & sweeps per step (converged)\\
			\midrule
			0.5 & 0.08 & 0.03 & 0.03 & 102\\
			1 & 0.61 & 0.04 & 0.03 & 55\\
			2 & 3.4 & 0.31 & 0.03 & 42\\
			5 & 14.4 & 4.3 & 0.03 & 88\\
			10 & 27.8 & 13.6 & 0.03 & 161\\
			50 & 60.0 & 46.9 & 0.03 & 691\\
			\bottomrule
		\end{tabular}
	\end{table}
	
	\subsection{P6: American put under Merton's jump-diffusion model}
\begin{figure}[t]
		\centering
		\includegraphics[width=\textwidth]{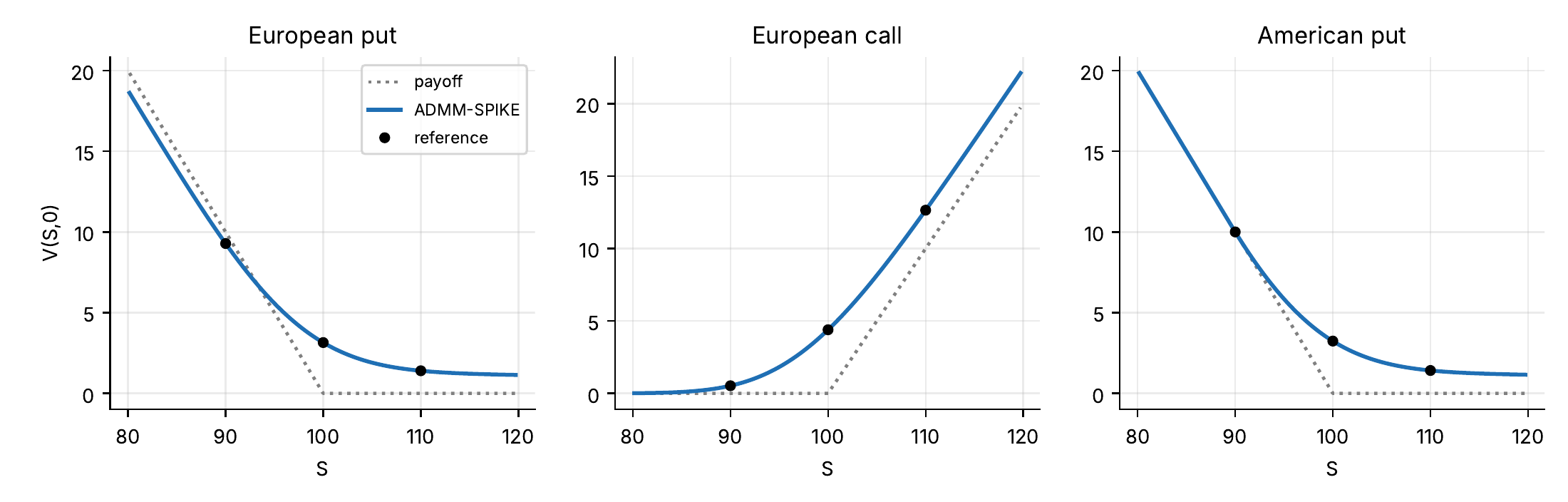}
		\caption{P6, P7: Merton jump-diffusion option values at $\tau=T=0.25$ against $S$ (ADMM-SPIKE, $M=800$, $N=400$), the payoff, and the reference values of \cite{maurya} at $S=90,100,110$.}
		\label{fig:fig6mertonprices}
	\end{figure}

\begin{figure}[t]
		\centering
		\includegraphics[width=\textwidth]{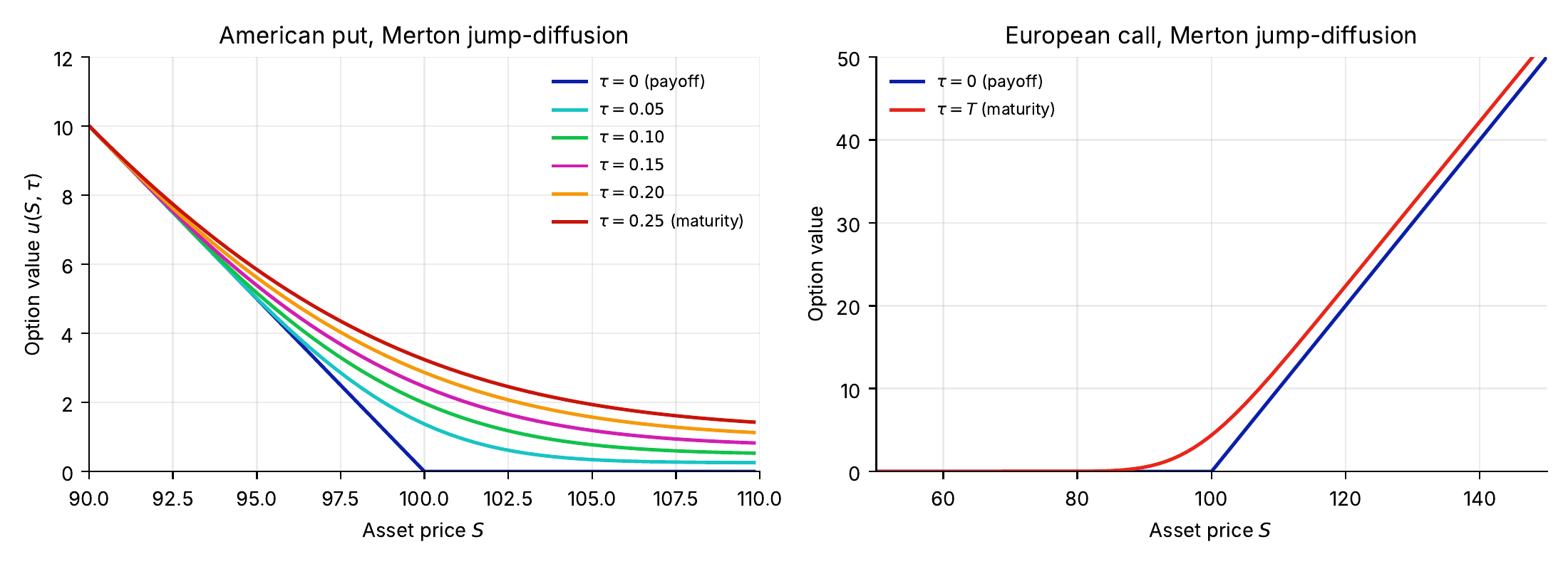}
		\caption{P6, P7: time evolution of the Merton jump-diffusion option values. Left: American put at $\tau=0,0.05,\dots,0.25$ (the curve at $\tau=0$ is the payoff); deep in the money the value stays at the exercise value $K-S$ at all times. Right: European call at $\tau=0$ (payoff) and at maturity $\tau=T=0.25$ over $50\le S\le150$. ADMM-SPIKE, $M=800$, $N=400$, $\rho=1$.}
		\label{fig:fig7mertontrajectories}
	\end{figure}

	With the settings of \cite{maurya} (domain $(-1.5,1.5)$, $M=800$ intervals, $N=400$ time steps, compensator $\kappa\approx-0.550$, $p=4$) the American put prices of Table~\ref{tab:merton} are within $2.6\times10^{-3}$ of the published values at $S=90,100,110$ (Figures~\ref{fig:fig6mertonprices} and~\ref{fig:fig7mertontrajectories}). The American-put ADMM needs about $40$ sweeps per step (with $\rho=1$, tolerance $10^{-10}$) and the solutions with $p=1$ and $p=4$ agree to $2.8\times10^{-14}$. The error of the American put falls with the grid from $2.2\times10^{-2}$ to $1.0\times10^{-3}$ ($M=200,\dots,1600$, $N=M/2$). The compensator $\lambda_J\kappa$ must be included in the drift; omitting it changes the model.
	
	\subsection{P7: European options and the tails of the integral operator}
	With the same settings the European prices of Table~\ref{tab:merton} are obtained. We recomputed the European prices independently with Merton's Poisson-mixture series \cite{merton}; they reproduce the published values to six decimals ($9.285418$, $3.149026$, $1.401186$ for the put and $0.527638$, $4.391246$, $12.643406$ for the call). The first-order time stepping gives errors of $1$--$2\times10^{-3}$ at $S=100$; the maximum error at $S=90,100,110$ falls with the grid ($M=200,400,800,1600$ with $N=M/2$) from $1.8\times10^{-2}$, $5.4\times10^{-3}$, $1.8\times10^{-3}$ to $6.5\times10^{-4}$ for the European put, with ratios $3.4$, $3.1$, $2.7$; the call behaves alike ($1.9\times10^{-2}\to7.0\times10^{-4}$). Our right-tail formula \eqref{eq:tailC} reproduces the call prices, which supports its use.
	\begin{table}[t]
		\centering\small
		\caption{P6, P7: Merton model, reference prices and the results of the present method ($M=800$, $N=400$).}\label{tab:merton}
		\begin{tabular}{lccc}
			\toprule
			Option & $S=90$ & $S=100$ & $S=110$\\
			\midrule
			European put, reference \cite{maurya} & 9.285418 & 3.149026 & 1.401186\\
			European put, present work & 9.286596 & 3.147264 & 1.401109\\
			European call, reference \cite{maurya} & 0.527638 & 4.391246 & 12.643406\\
			European call, present work & 0.528727 & 4.389396 & 12.643240\\
			American put, reference \cite{maurya} & 10.003822 & 3.241251 & 1.419803\\
			American put, present work & 10.003903 & 3.238714 & 1.419601\\
			\bottomrule
		\end{tabular}
	\end{table}
	
	\subsection{Summary of the verification}
	In every case where an independent reference exists, the converged ADMM-SPIKE iteration agrees with it to the discretization error (exact solutions in 1D, binomial and closed-form prices, published Merton prices), and the SPIKE solution does not depend on the number of processes to $3.5\times10^{-13}$ or better. The experiments also show that the inner iteration of the option problems must be run to convergence, or that the augmentation parameter must be of the size of the diagonal of the time-stepping matrix when the number of sweeps is fixed (Table~\ref{tab:rho}).
	
	%======================================================================
	\section{Summary and conclusions}\label{sec:concl}
	In this study, we have presented split Bregman (ADMM) solvers for obstacle problems and for the LCPs that arise in the pricing of American options, in which all linear algebra is a SPIKE-based parallel banded solve. The factorization, the spikes and the reduced system are computed once and reused in every iteration and every time step, so that only a pointwise projection and a multiplier update remain. The method was implemented in Fortran~90 with MPI. A Lyapunov functional of the primal--dual pair shows that the iteration converges for every penalty parameter $\lambda>0$, for the symmetric obstacle problem and for the non-symmetric strongly monotone LCPs of option pricing alike, and SPIKE leaves the iterates unchanged for every number of processes. In every case with an independent reference, the converged iteration agrees with it to the discretization error: with the exact solutions of the one-dimensional obstacle problems, with the closed-form Black--Scholes call (error $\le1.5\times10^{-3}$), with the binomial American put ($V(100)=6.0885$ against $6.0902$), and with the published Merton prices for the American put and for European puts and calls (within $2.6\times10^{-3}$, decreasing with the grid).
	
	Compared with PSOR, penalty and active-set methods, the method has the following advantages. One fixed matrix is factorized once, and since the contact set does not enter the matrix, nothing is refactorized. No relaxation or penalty parameter has to be tuned for convergence. The projection and the multiplier update are pointwise and perfectly parallel, while SPIKE distributes the linear solve. A single algorithm and a single convergence proof cover symmetric obstacle problems and non-symmetric option LCPs under the Black--Scholes and Merton jump-diffusion models, and the results are independent of the number of processes ($3.5\times10^{-13}$ or better for $1,2,4,8$ processes). The two-dimensional solve scales almost linearly up to four processes; for $n_x=256$ the time falls from $1343.5$~s to $343.1$~s, a speed-up of $3.92$.
	
	The limitations of the method are the following. The inner iteration has to be run to convergence in every time step, or the augmentation parameter has to be of the size of the diagonal of the time-stepping matrix: with $\rho=50$ and twelve sweeps per step the American put at the money is $3.23$ instead of $6.09$ ($47\,\%$ error), whereas $\rho=1$ with the same sweeps gives $0.04\,\%$. ADMM needs about $10^3$ iterations, whereas active-set methods need only a few, and no convergence rate is proved. For the two-dimensional obstacle problem the stopping rule $\|u^{k+1}-u^k\|_\infty<5\times10^{-4}$ is not an error estimate, since it is met while the error is still $0.18$. The parallel speed-up saturates beyond four processes, and in one dimension the diagonal blocks were treated as dense matrices, so that the timings do not support a speed-up claim; for a tridiagonal system SPIKE cannot be expected to beat a serial Thomas solve. Future work includes a serial banded-LU baseline on the same machine, a residual-based stopping test in two dimensions, banded block factorizations in one dimension, an investigation of the saturation beyond four processes, and the numerical treatment of Kou's double-exponential jump model, for which the tail formula is already given.
	
	%======================================================================
	
\end{document}